\documentclass[10pt]{article}
\usepackage{amssymb,amsmath,amsthm}
\usepackage{color}
\newcommand{\R}{\mathbb{R}}
 \UseRawInputEncoding

\begin{document}
\title{\Large\bf{ Existence and convergence of ground state solutions for a $(p,q)$-Laplacian system on weighted graphs}}
\date{}
\author {  Xuechen Zhang,$^{1}$       Xingyong Zhang$^{1,2}$\footnote{Corresponding author, E-mail address: zhangxingyong1@163.com}\\
       {\footnotesize $^1$Faculty of Science, Kunming University of Science and Technology,}\\
       {\footnotesize Kunming, Yunnan, 650500, P.R. China.}\\
       {\footnotesize $^2$Research Center for Mathematics and Interdisciplinary Sciences, Kunming University of Science and Technology,} \\
       {\footnotesize Kunming, Yunnan, 650500, P.R. China.}}

 \date{}
 \maketitle

 \begin{center}
 \begin{minipage}{15cm}
 \par
 \small  {\bf Abstract:} We investigate the existence of ground state solutions for a $(p,q)$-Laplacian system with $p,q>1$ and potential wells
on a weighted locally finite graph $G=(V,E)$. By making use of the method of Nehari manifold and the Lagrange multiplier rule, we prove that if the nonlinear term $F$ takes on the super-$(p, q)$-linear growth and the potential functions $a(x)$ and $b(x)$ satisfy some suitable conditions, then for any fixed parameter $\lambda\geq1$, the system is provided with a ground state solution $(u_\lambda, v_\lambda)$. Additionally, we set up the convergence property of the solutions set $\{(u_\lambda, v_\lambda)\}$  when  $\lambda \rightarrow +\infty$.

 \end{minipage}
 \end{center}
  \allowdisplaybreaks
 \vskip2mm
 {\section{Introduction }}
\setcounter{equation}{0}
 In recent years, the analysis on graphs has attracted some attentions.  For example, in \cite{Bauer F 2015}, Bauer-Horn-Lin-Lippner-Mangoubi-Yau proved a discrete version of Li-Yau inequality valid for solutions to the heat equation on graphs. In \cite{Horn P 2019},  Horn-Lin-Liu-Yau proved a Gaussian estimate for the heat kernel, along with Poincar\'e and Harnack inequalities. It is remarkable that some analysis on graphs has been applied to the investigation of machine learning, data analysis, neural network, image processing etc (for example, see \cite{Alkama S2014}, \cite{Arnaboldi V 2015}, \cite{Elmoataz A2015}, \cite{Ta V T 2010}, \cite{Ta V T 2008}).

\par
Next, we would like to address some investigations of the partial differential equations on graphs. For example, the heat equations \cite{Huang X 2012}, the Liouville type equations \cite{Ge H 2018}, the Fokker-Planck equation \cite{Chow S N 2017} and the Schr\"{o}dinger equation \cite{Chow S N 2019}. Especially, in \cite{Grigor'yan A 2017}, Grigor-Lin-Yang studied the following nonlinear equation on a locally finite graph $G=(V,E)$:
\begin{eqnarray*}
\label{F1}
-\Delta u+hu=f(x,u) \;\;\;\;\hfill \mbox{in}\;\; V,
\end{eqnarray*}
where $\Delta$ is a graph Laplacian defined (\ref{Eq4}) below, $V$ is the  vertex set and $E$ is the edge set,  $h:V\to \R$ satisfies some suitable assumptions and $f$ satisfies the superquadratic growth condition. They established some existence results via mountain pass theorem.  In \cite{Grigor'yan A 2016 yamabe}, Grigor-Lin-Yang established some existence results to the Yamabe type equation
\begin{eqnarray*}
\label{F2}
 \begin{cases}
  -\Delta u- \alpha u= |u|^{p-2}u,\;\;\;\hfill \mbox{in $\Omega^{\circ}$},\\
  u=0 ,\;\;\;\;\;\;\;\;\;\;\;\;\;\;\;\;\;\;\;\;\;\;\;\;\hfill \mbox{on $\partial\Omega$ },
   \end{cases}
\end{eqnarray*}
by the mountain pass theorem, where $p \geq 2, \alpha < \lambda_1 (\Omega)$ which is the first eigenvalue of $-\Delta$. They also considered a $p$-Laplacian equation with Dirichlet boundary value and a generalized poly-Laplacian equation with Dirichlet boundary value on weighted locally finite graph, and a $p$-Laplacian equation  and a generalized poly-Laplacian equation on finite graph by the same method.
\par
Besides, the existence of ground state solutions on graphs have also attracted some academics via the Nehari manifold method which was developed by Nehari firstly in \cite{Nehari Z 1959}. In \cite{Zhang N 2018}, Zhang-Zhao studied the convergence of ground state solutions for the following nonlinear Schr\"{o}dinger equation on a locally finite graph $G=(V,E)$,
\begin{eqnarray*}
\label{F4}
-\Delta u+ (\lambda a+1) u= |u|^{p-1}u \;\;\;\;\hfill \mbox{in V},
\end{eqnarray*}
 where $\lambda >1$ and $a(x)$ satisfies the following assumptions:\\
$(a_1) \;\;a(x)\geq 0$ and the potential well $\Omega= \{x \in V : a(x)=0\}$ is a non-empty, connected and bounded
domain in $V$;\\
$(a_2)$ \;\;there exists a vertex $x_0\in V$ such that $a(x) \rightarrow 0$ as $d(x,x_0)\rightarrow +\infty$.
\vskip2mm
\noindent
And, as $\lambda \rightarrow +\infty$, the solutions family $\{u_\lambda\}$ converges to a ground state solution of the following Dirichlet problem
\begin{eqnarray*}
\label{F5}
\begin{cases}
  -\Delta u+ u= |u|^{p-1}u\;\;\;\;\hfill \mbox{in}\;\; \Omega,\\
  u=0 ,\;\;\;\;\hfill \mbox{on}\;\; \partial\Omega .
   \end{cases}
\end{eqnarray*}
In \cite{Han X 2020}, Han-Shao-Zhao considered the following biharmonic equation on locally locally finite graph $G=(V,E)$:
\begin{eqnarray*}
\label{F6}
\Delta^2 u-\Delta u+ (\lambda a+1) u= |u|^{p- 2}u \;\;\;\;\hfill \mbox{in}\;\; V,
\end{eqnarray*}
where $\lambda >1$, $a(x)$ satisfied $(a_1)$ and $(a_2)$ and they investigated the existence of ground state solution $u_\lambda$ on graphs by the Nehari manifold method. Moreover, as $\lambda \rightarrow +\infty$,  they also obtained that the solution $u_\lambda$ converges to a ground state solution of the  following Dirichlet problem
\begin{eqnarray*}
\label{F7}
\begin{cases}
  \Delta^2 u-\Delta u+ u= |u|^{p-1}u\;\;\;\;\hfill \mbox{in}\;\; \Omega,\\
  u=0 ,\;\;\;\;\hfill \mbox{on}\;\; \partial\Omega .
   \end{cases}
\end{eqnarray*}
In \cite{Han X L 2021}, Han-Shao investigated the following $p$-Laplacian equation with $p\geq2$ on locally finite graph $G=(V,E)$:
\begin{eqnarray*}
\label{F8}
-\Delta_p u+ (\lambda a+1)u^{p-2}u= f(x,u)\;\;\;\;\hfill \mbox{in}\;\; V,
\end{eqnarray*}
where $\lambda >1$. They got that the equation admits a ground state solution $u_\lambda$ via the Nehari manifold method and the deformation lemma. Besides, they also proved the convergence of ground state solution. They assumed that the nonlinear term $f$ satisfies the following conditions:\\
$(f_1)$ \;\;for any $x\in V$, $f(x,s)$ is continuous in $s \in \R$, $f(x,0)=0$, and for any fixed $M>0$ there exists a constant $A_M$ such that $\max_{|s|<M}f(x,s)\leq A_M$ for all $x \in V$;\\
$(f_2)$ \;\;there exists some $\alpha > p$ such that for any $s \in \R\setminus\{0\}$ and $x \in V$ there holds
$$
0< \alpha F(x,s):=\alpha\int^s_0f(x,t)dt\leq sf(x,s);
$$
$(f_3)$ \;\;for any $x \in V$, there holds
$$
\limsup\limits_{|s|\rightarrow 0}\frac{|f(x,s)|}{|s|^{p-1}} < \lambda_p :=\inf\limits_{u\not\equiv 0}\frac{\int_V(|\nabla u|^p+(\lambda a+1)|u|^p)d\mu}{\int_V|u|^pd\mu};
$$
$(f_4)$ \;\;there exist some $q>p$ and $C>0$ such that
\begin{eqnarray*}
|f(x,s)| \leq C(1+|s|^{q-1})\;\;\;\; \mbox{uniformly in}  \;\;x\in V;
\end{eqnarray*}
$(f_5)$ \;\;$s \mapsto \frac{f(x,s)}{|s|^{p-1}}$ is strictly increasing on $(- \infty, 0)$ and $(0, +\infty)$ for all $x\in V$.\\
The assumptions on the potential $a(x)$ satisfies $(a_1)$ and\\
$(a_2)'$ \;\;$(a(x)+1)^{-1}\in L^{\frac{1}{p-1}}(V)$.

\vskip2mm
\noindent
In \cite{Shao 2023}, Shao studied the following $p$-Laplacian systems with $p\geq2$ on locally finite graph $G=(V,E)$:
\begin{eqnarray*}
\label{f1}
 \begin{cases}
  -\Delta_p u+(\lambda a+1)|u|^{p-2}u=\frac{1}{\kappa}F_u(u,v),\;\;\;\;\hfill \mbox{in}\;\; V,\\
  -\Delta_p v+(\lambda b+1)|v|^{p-2}v=\frac{1}{\kappa}F_v(u,v),\;\;\;\;\hfill \mbox{in}\;\;V,\\
   \end{cases}
\end{eqnarray*}
where $\lambda>1, 2\leq p<\kappa< \infty$ and $F$ satisfies the following conditions:\\
$(F)$ \;\;$F \in C^1(\R^2,\R^+)$ and $F(tu,tv)=t^\kappa F(u,v)$ holds for all $(u,v)\in \R^2$, where $t>0$.\\
The assumptions on the potential $a(x)$ and $b(x)$ satisfies\\
$(a_3)$ \;\;$a(x)\geq 0, b(x) \geq 0$. The potential well $\Omega_a:=\{x\in V:a(x)=0\}$, $\Omega_b:=\{x\in V:b(x)=0\}$, and $\Omega_a, \Omega_b$ and $\Omega_a\cap\Omega_b$ are all non-empty, connected and bounded domains in $V$.\\
$(a_4)$ \;\;there exists a vertex $x_0\in V$ such that $a(x) \rightarrow 0$ and $b(x) \rightarrow 0$  as $d(x,x_0)\rightarrow +\infty$.\\
Shao established the existence of ground state solutions by the Nehari manifold method. When $\lambda \rightarrow +\infty$, the ground state solutions family converges to ground state solutions of the corresponding Dirichlet equation.

\par
The ground state solutions for partial differential equations in the Euclidean setting have been studied extensively. We mainly refer to the work in \cite{Papageorgiou N S 2020} which inspired our work.  In \cite{Papageorgiou N S 2020}, Papageorgiou-R\u{a}dulescu established the existence result of the ground state solution by the Nehari manifold method and the Lagrange multiplier rule for the following double phase problems:
\begin{eqnarray*}
\label{e15}
 \begin{cases}
  -\Delta_p u-\mbox{div}(a(z)|Du|^{q-2}Du)=f(x,u),\;\;\;\;\hfill \mbox{in}\;\; \Omega,\\
  u=0,\;\;\;\;\hfill \mbox{in}\;\;\partial\Omega, 1<q<p,\\
   \end{cases}
\end{eqnarray*}
where $\Omega$ is a bounded domain in $\R^N(N\geq 2)$ with a smooth boundary, $p^*$ is the critical Sobolev exponent corresponding to $p$, which is defined by
\begin{eqnarray*}
p^*=
 \begin{cases}
  \frac{Np}{N-p},& \text {if $p<N$}\\
  +\infty,&  \text {if $N\leq p $,}
   \end{cases}
\end{eqnarray*}
$f$ satisfies the following assumptions:\\
$H(f)$ \; $f: \Omega \times \R\rightarrow \R$ is a measurable function such that for a.e. $z\in \Omega, f(z,0)=0, f(z,\cdot)\in C^1(\R\backslash\{0\})$ and       \\
$(i)$ \; $|f'_x(z,x)|\leq a_0(z)(1+|x|^{r-2})$ for a.e. $z\in \Omega$ and all $x\in \R$, with $a_0\in L^\infty(\Omega), 1<q<p<r<p^*;$\\
$(ii)$ \;Let $F(z,x)= \int_0^xf(z,x)ds$. Then $\lim\limits_{x\pm\infty}\frac{F(z,x)}{|x|^p}=+\infty$ uniformly for a.e. $z\in \Omega$ and there exist constants $\tau\in ((r-p)\max\{1,\frac{N}{p}\} , p^*)$ and $\beta_0>0$ such that
$$
\beta_0\leq \liminf\limits_{x\rightarrow \pm\infty}\frac{f(z,x )x -p F(z,x)}{|x|^{\tau_1}}
\; \mbox{ for all} \;\; (x, t,s )\in V\times \R^2;
$$
$(iii)$ \; $\lim\limits_{x\rightarrow 0}\frac{f(z,x)}{ |x|^{q-2}x}=0$ for all $x\in \Omega$;\\
$(iv)$ \; $0< (p-1)f(z,x)x \leq f'_x(z,x)x^2  $ for a.e. $x\in \Omega$ and all $x\neq 0$.

\vskip2mm
\par
Motivated by \cite{Han X L 2021}, \cite{Grigor'yan A 2016 yamabe}, \cite{Papageorgiou N S 2020} and \cite{Shao 2023}, we consider the following $(p,q)$-Laplacian systems
\begin{eqnarray}
\label{eq1}
 \begin{cases}
  -\Delta_p u+(\lambda a+1)|u|^{p-2}u=F_u(x,u,v),\;\;\;\;\hfill \mbox{in}\;\; V,\\
  -\Delta_q v+(\lambda b+1)|v|^{q-2}v=F_v(x,u,v),\;\;\;\;\hfill \mbox{in}\;\;V,\\
   \end{cases}
\end{eqnarray}
where $V$ is a locally finite graph, $p,q>1$, $\lambda\geq1$,  and $a(x)$ and $b(x)$ are potential functions. The nonlinear term $F:V\times \R\times \R\to \R$ and $\Delta_s, s=p,q$ is the discrete $s$-Laplacian on graphs. We prove that for every given $\lambda>1$,  the systems admits a ground state solution $(u_\lambda, v_\lambda)$ via the method of Nehari manifold and the Lagrange multiplier rule. As $\lambda \rightarrow +\infty$, we also obtained that the solutions family $\{(u_\lambda, v_\lambda)\} $ converges to the solution of the limit problem defined on the potential wells $\Omega_a$ and $\Omega_b$:
\begin{eqnarray}
\label{e5}
 \begin{cases}
 -\Delta_p u+|u|^{p-2}u=F_u(x,u,v)\;\;\;\;\hfill \mbox{in}  \;\; \Omega_a, \\
  -\Delta_q v+|v|^{q-2}v=F_v(x,u,v)\;\;\;\;\hfill \mbox{in}  \;\; \Omega_b, \\
  u=0\;\;\;\;\hfill \mbox{on}  \;\; \partial\Omega_a, \\
  v=0\;\;\;\;\hfill \mbox{on}  \;\; \partial\Omega_b,
   \end{cases}
\end{eqnarray}
where $\Omega_a=\{x \in V: a(x)=0\}$ and $\Omega_b=\{x \in V: b(x)=0\}$, $\partial \Omega_a$ denotes the boundary of $\Omega_a$ and $\partial \Omega_b$ denotes  the boundary of $\Omega_b$ which  is defined by (\ref{e16}) below.

\par
To describe our problems and results more clearly, we review  some concepts and assumptions (see \cite{Han X L 2021}, \cite{Grigor'yan A 2016 yamabe}, \cite{Grigor'yan A 2017}). Let $G=(V, E)$ be a graph, where $V$ denotes the vertex set and $E$ denotes the edge set. For any $x \in V$, if there are only finite $y \in V$ such that $xy \in E$, then $G$ is called a locally finite graph. $G$ is connected if any two vertices $x$ and $y$ can be connected via finitely many edges. Let a measure $\mu:V\rightarrow \R^+$ be finite. If there exists a constant $\mu_{\min}>0$ such that $\mu(x)\geq \mu_{\min}$ for all $x\in V$, we call that $\mu$ is a uniformly positive measure. $G$ is said to be symmetric, namely $\omega_{xy}=\omega_{yx}$. For any edge $xy\in E$ with two vertexes of $x,y\in V$, assume that its weight $\omega_{xy}>0$ and for any $x\in V$, $\sum_{y\thicksim x}\omega_{xy} <C$, where $C$ is a positive constant. Here and throughout this paper, $y \thicksim  x$ stands for any vertex $y$ connected with $x$ by an edge $xy \in E$.
The distance $d(x,y)$ of two vertices $x,y \in V$ is defined by the minimal number of edges which connect these two vertices. $\Omega \subset V$ is a bounded domain in $V$, if the distance $d(x,y)$ is uniformly bounded from above for any $x,y \in \Omega$. Obviously, a bounded domain of a locally finite graph contains only finite vertices. Denote that the boundary of $\Omega$ by
\begin{eqnarray}
\label{e16}
\partial \Omega :=\{y \in V\setminus\Omega : x \in \Omega \;\;\mbox{such that}\;\; xy \in  E\}
\end{eqnarray}
and the interior of $\Omega$ by $\Omega^\circ$. Obviously, $\Omega^\circ=\Omega$.
\par
For any $x\in V$, we define
\begin{eqnarray*}
\label{Eq4}
\Delta \psi(x)=\frac{1}{\mu(x)}\sum\limits_{y\thicksim x}w_{xy}(\psi(y)-\psi(x)),
\end{eqnarray*}
where $\psi : V\rightarrow \R$. The corresponding gradient form is
\begin{eqnarray*}
\label{Eq5}
\Gamma(\psi_1,\psi_2)(x)=\frac{1}{2\mu(x)}\sum\limits_{y\thicksim x}w_{xy}(\psi_1(y)-\psi_1(x))(\psi_2(y)-\psi_2(x)).
\end{eqnarray*}
Write $\Gamma(\psi)=\Gamma(\psi,\psi)$. The length of the gradient is defined by
\begin{eqnarray*}
\label{Eq6}
|\nabla \psi|(x)=\sqrt{\Gamma(\psi)(x)}=\left(\frac{1}{2\mu(x)}\sum\limits_{y\thicksim x}w_{xy}(\psi(y)-\psi(x))^2\right)^{\frac{1}{2}}.
\end{eqnarray*}
For any function $\psi:V\rightarrow\mathbb{R}$, we denote
\begin{eqnarray*}\label{Eq8}
\int\limits_V \psi(x) d\mu=\sum\limits_{x\in V}\mu(x)\psi(x).
\end{eqnarray*}
When $p\geq2$, we define the $p$-Laplacian operator by $\Delta_p\psi$ with
\begin{eqnarray*}\label{Eq9}
\Delta_p\psi(x)=\frac{1}{2\mu(x)}\sum\limits_{y\sim x}\left(|\nabla \psi|^{p-2}(y)+|\nabla \psi|^{p-2}(x)\right)\omega_{xy}(\psi(y)-\psi(x)).
\end{eqnarray*}
In the distributional sense, $\Delta_p \psi$ can be written as follows. For any $\phi\in\mathcal{C}_c(V)$,
\begin{eqnarray*}\label{Eq10}
\int\limits_V(\Delta_p \psi)\phi d\mu=-\int\limits_V|\nabla \psi|^{p-2}\Gamma(\psi,\phi)d\mu,
\end{eqnarray*}
where $\mathcal{C}_c(V):=\{u:V\rightarrow \R:\{x \in V: u(x)\neq 0\} \mbox{ is of finite cardinality}\}$.

\par
When $1\leq\theta <+\infty$, we define
$$
L^{\theta }(V)=\left\{\psi:V\to\R\Big|\int_{V}|\psi(x)|^{\theta } d\mu<\infty\right\}
$$
and the norm  by
\begin{eqnarray*}
\label{b1}
\|\psi\|_{{\theta },V}=\left(\int_{V}|\psi(x)|^{\theta } d\mu\right)^\frac{1}{{\theta }}.
\end{eqnarray*}
When ${\theta } =+\infty$, we define
$$
L^\infty(V)=\left\{\psi:V\to\R\Big|\sup_{x\in V}|\psi(x)|<\infty\right\}
$$
with the norm
\begin{eqnarray*}
\label{Z7}
\|\psi\|_{\infty,V}=\sup_{x \in V}|\psi(x)|< \infty.
\end{eqnarray*}
We define $\|(\psi,\phi)\|_{L(V)}= \|\psi\|_{\theta_1}+\|\phi\|_{\theta_2}$, where $1\leq \theta_i \leq +\infty, i=1,2.$

\vskip2mm
\par
 Next, we introduce the following assumptions on $F:V\times\R\times\R\rightarrow\R$ and  \\
$(F_1)$ \;  $F(x,0,0)=0$ and $F(x,t,s)$ is twice continuously differentiable in $(t,s)\in \R^2$ for all $x\in V$;\\
$(F_2)$ \; there exist constants $\alpha>\beta:=\max\{p,q\}$ and $0  <\varrho < \min\{\frac{\alpha-p}{p }, \frac{\alpha-q}{q }\}$ such that
$$
\alpha F(x,t,s)\leq F_t(x,t,s)t+F_s(x,t,s)s+\varrho ( |t|^p+|s|^q)
\; \mbox{ for all} \;\; (x, t,s )\in V\times \R^2,
$$
where $F_t(x,t,s)=\frac{\partial F(x,t,s)}{\partial t}$ and  $F_s(x,t,s)=\frac{\partial F(x,t,s)}{\partial s}$;\\
$(F_3)$ \; there exist constants $r_1 , r_2>\alpha$, two functions $C_i:V\rightarrow \R, i=1,2$ such that $C_1\in L^1(V)\cap L^\infty(V)$ with $0<C_1(x)\leq\frac{1}{1+\beta}$ for all $x\in V$, $C_2\in L^1(V)\cap L^\infty(V)$ with $C_2(x)>0$ for all $x\in V$ and
$$
|F_t(x,t,s)|\leq C_1(x)(|t|^{p-1}+|s|^{\frac{q(p-1)}{p}})+C_2(x)(|t|^{r_1-1}+|s|^{\frac{r_2(r_1-1)}{r_1}}),   \; \mbox{for all $x\in V$ and $\forall(s,t)\in \R^2$};
$$
$$
|F_s(x,t,s)|\leq C_1(x)(|t|^{\frac{p(q-1)}{q}}+|s|^{q-1})+C_2(x)(|t|^{\frac{r_1(r_2-1)}{r_2}}+|s|^{r_2-1}),  \; \mbox{for all $x\in V$ and $\forall(s,t)\in \R^2$};
$$
$(F_4)$ \; $F_{ts}(x,t,s)ts\geq0$ and $0< (\beta-1)(F_t(x,t,s)t+ F_s(x,t,s)s)< F_{tt}(x,t,s)t^2+ F_{ss}(x,t,s)s^2$ for all $x\in V$ and $\forall (t,s) \in \R^2\setminus\{(0,0)\}$.

\vskip2mm
\par
The assumptions on the potential functions $a(x)$ and  $b(x)$ are:\\
$(A_1)$ \; $a(x)\geq 0, b(x) \geq 0$ for all $x\in V$. The potential wells $\Omega_a:=\{x\in V:a(x)=0\}$, $\Omega_b:=\{x\in V:b(x)=0\}$,  and $\Omega_a, \Omega_b$ and $\Omega_a\cap\Omega_b$ are all non-empty, connected and bounded domains in $V$.\\
$(A_2)$ \; $(a(x)+1)^{-1}\in L^{\frac{1}{p-1}}(V)$ and $(b(x)+1)^{-1}\in L^{\frac{1}{q-1}}(V)$.

\vskip2mm
\par
We define
\begin{eqnarray*}
W^{1,s}(V)=\left\{u:V\to\R\Big|\int_V(|\nabla  u|^s+|u|^s)d\mu<\infty\right\}
\end{eqnarray*}
which is provided with the norm
\begin{eqnarray*}
\label{e1}
\|u\|_{W^{1,s}(V)}=\left(\int_V(|\nabla  u|^s+|u|^s)d\mu\right)^\frac{1}{s},
\end{eqnarray*}
where $s>1$. Let us consider the space $W:=W^{1,p}(V)\times W^{1,q}(V)$ and  $\|(u,v)\|_W=\|u\|_{W^{1,p}(V)}+\|v\|_{W^{1,q}(V)}$.
\par
Define
$$
W_\lambda(a):=\left\{u\in W^{1,p}(V)\Big|\int_V\lambda a|u|^pd\mu<\infty\right\}
$$
and
$$
W_\lambda(b):=\left\{v\in W^{1,q}(V)\Big|\int_V\lambda b|v|^qd\mu<\infty\right\}.
$$
\par
To study the problem (\ref{eq1}), it is natural to consider the function space
$$
W_\lambda:=\left\{(u,v)\in W\Big|\int_V(\lambda a|u|^p+\lambda b|v|^q)d\mu<\infty\right\}.
$$
It is easy to see that the space $W_\lambda:=W_\lambda(a)\times W_\lambda(b)$. Define the norm $\|(u,v)\|_{W_\lambda}=\|u\|_{W_\lambda(a)}+\|v\|_{W_\lambda(b)}$, where
$$
\|u\|_{W_\lambda(a)}=\left(\int_V(|\nabla u|^p+(\lambda a+1)|u|^p)d\mu\right)^{\frac{1}{p}}
$$
and
$$
\|v\|_{W_\lambda(b)}=\left(\int_V(|\nabla v|^q+(\lambda b+1)|v|^q)d\mu\right)^{\frac{1}{q}}.
$$
Then, $(W_\lambda,\|\cdot\|_{W_\lambda})$ is a reflexive Banach space (see \cite{Shaomeng 2023}).
\par
The functional related to (\ref{eq1}) is defined by
\begin{eqnarray*}
\label{e2}
J_\lambda(u,v)=\frac{1}{p}\int_V(|\nabla u|^p+(\lambda a+1)|u|^p)d\mu+\frac{1}{q}\int_V(|\nabla v|^q+(\lambda b+1)|v|^q)d\mu-\int_V F(x,u,v)d\mu.
\end{eqnarray*}
Under the assumptions $(F_1)$ and $(F_3)$, a standard procedure can show that $J_\lambda\in C^1(W_\lambda,\R)$ and
\begin{eqnarray*}
\label{e3}
\langle J'_\lambda(u,v),(\phi_1,\phi_2)\rangle
&=&\int_V(|\nabla u|^{p-2} \Gamma(u,\phi_1)+(\lambda a+1)|u|^{p-2}u\phi_1)d\mu
+\int_V(|\nabla v|^{q-2} \Gamma(v,\phi_2)+(\lambda b+1)|v|^{q-2}v\phi_2)d\mu\nonumber\\
&&-\int_V (F_u(x,u,v)\phi_1+F_v(x,u,v)\phi_2)d\mu
\end{eqnarray*}
for all $(u,v), (\phi_1, \phi_2)\in W_\lambda$.
By Lemma 3.1 below, the critical point of $J$ is the point-wise solution of system (\ref{eq1}).
\par
Define the Nehari manifold
$$
\mathcal{N_\lambda}:=\{(u,v)\in W_\lambda\backslash\{(0,0)\}:\langle J'_\lambda(u,v),(u,v)\rangle=0\}
$$
and the least energy level
\begin{eqnarray*}
\label{e4}
m_\lambda:=\inf_{(u,v)\in \mathcal{N_\lambda}}J_\lambda(u,v).
\end{eqnarray*}
\par
It is well known that if $m_\lambda$ can be achieved by some function $(u_\lambda, v_\lambda)\in \mathcal{N_\lambda}$ and $(u_\lambda, v_\lambda)$ is a critical point of the functional $J_\lambda$, then $(u_\lambda, v_\lambda)$ is  a ground state solution of (\ref{eq1}). Next, we provide our main results.

\vskip2mm
\noindent
{\bf Theorem 1.1.} {\it Let $G=(V, E)$ be a locally finite graph and it is connected, symmetric and there exists $\mu_{\min}>0$ such that $\mu(x)\geq \mu_{\min}$ for all $x\in V$. For any $x\in V$, $\sum_{y\thicksim x}\omega_{xy} <C$, where $C$ is a positive constant. Assume $\frac{1}{\gamma}<\beta +\frac{1}{\beta} -1$, where $\gamma=\min\{p,q\}$,  and $(A_1), (A_2)$ and $(F_1)-(F_4)$ hold. Then for every $\lambda\geq1$, there exists a ground state solution $(u_\lambda, v_\lambda)$ of system (\ref{eq1}) and $\|(u_\lambda, v_\lambda)\|_{W_\lambda}$ satisfies (\ref{**}) below.}

\vskip2mm
\noindent
{\bf Remark 1.1.} There exists examples satisfying Theorem 1.1. For example, let $\alpha>\beta:=\max\{p,q\}$ and
\begin{eqnarray*}
&&C_1(x)=\frac{1}{(1+\beta)(1+|x|)^2},\;\;\;\;\;\;C_2(x)=\frac{2 \beta}{(1+|x|)^2},\\
&&F(x,t,s)=\frac{|t|^{\alpha}+|s|^{\alpha}}{2\alpha(1+\beta)(1+|x|)^2}\left(1-\frac{1}{\ln(e^2+|t|^\alpha+|s|^\alpha)}\right).
\end{eqnarray*}
\par

\par
In order to investigate  the asymptotic behavior of $(u_\lambda, v_\lambda)$ to the solution of (\ref{e5}), we define $W_{\Omega_a}$ by the completion of $C_c(\Omega_a)$ under the norm
$$
\|u\|_{W_{\Omega_a}}=\left(\int_{\Omega_a\cup\partial\Omega_a}|\nabla u|^pd\mu+
\int_{\Omega_a}|u|^pd\mu\right)^{\frac{1}{p}},
$$
and $W_{\Omega_b}$ by the completion of $C_c(\Omega_b)$ under the norm
$$
\|v\|_{W_{\Omega_b}}=\left(\int_{\Omega_b\cup\partial\Omega_b}|\nabla v|^qd\mu+
\int_{\Omega_b}|v|^qd\mu\right)^{\frac{1}{q}},
$$
where $C_c(\Omega)$ denotes the set of all functions $u:\Omega\rightarrow \R$ satisfying $supp\;u\subset \Omega$ and $u=0$ on $V \setminus A$ for some subset $\Omega \subset V$. It is suitable to study (\ref{e5}) in the space $W_\Omega:=W_{\Omega_a}\times W_{\Omega_b}$ and $W_\Omega$ is a reflexive Banach space with finite dimensional. Define the norm $\|(u,v)\|_{W_\Omega}=\|u\|_{W_{\Omega_a}}+\|v\|_{W_{\Omega_b}}$.
\par
We define when $1\leq\theta_1<+\infty$,
$$
L^{\theta_1}(\Omega_a)=\left\{\psi:\Omega_a\to\R\Big|\int_{\Omega_a}|\psi(x)|^{\theta_1} d\mu<\infty\right\}
$$
and the norm  by
\begin{eqnarray*}
\label{bb1}
\|\psi\|_{{\theta_1},\Omega_a}=\left(\int_{\Omega_a}|\psi(x)|^{\theta_1} d\mu\right)^\frac{1}{\theta_1}.
\end{eqnarray*}
When ${\theta_1} =+\infty$, we define
$$
L^\infty(\Omega_a)=\left\{\psi:\Omega_a\to\R\Big|\sup_{x\in \Omega_a}|\psi(x)|<\infty\right\}
$$
with the norm
\begin{eqnarray*}
\label{Zz7}
\|\psi\|_{\infty,\Omega_a}=\sup_{x \in \Omega_a}|\psi(x)|< \infty.
\end{eqnarray*}
Similarly, when $1\leq\theta_2<+\infty$,
$$
L^{\theta_2}(\Omega_b)=\left\{\phi:\Omega_b\to\R\Big|\int_{\Omega_b}|\phi(x)|^{\theta_2} d\mu<\infty\right\}
$$
and the norm  by
\begin{eqnarray*}
\label{bb1}
\|\phi\|_{{\theta_2},\Omega_b}=\left(\int_{\Omega_b}|\phi(x)|^{\theta_2} d\mu\right)^\frac{1}{\theta_2}.
\end{eqnarray*}
When ${\theta_2} =+\infty$, we define
$$
L^\infty(\Omega_b)=\left\{\phi:\Omega_b\to\R\Big|\sup_{x\in \Omega_b}|\phi(x)|<\infty\right\}
$$
with the norm
\begin{eqnarray*}
\label{Zz7}
\|\phi\|_{\infty,\Omega_b}=\sup_{x \in \Omega_b}|\phi(x)|< \infty.
\end{eqnarray*}
We also define $\|(\psi,\phi)\|_{\theta,\Omega}= \|\psi\|_{\theta_1,\Omega_a}+\|\phi\|_{\theta_2,\Omega_b}$.
\par
The functional related to (\ref{e5}) is
\begin{eqnarray*}
\label{e6}
J_\Omega(u,v)=\frac{1}{p}\left(\int_{\overline{\Omega}_a} |\nabla u|^pd\mu+\int_{\Omega_a}|u|^p d\mu\right)+\frac{1}{q}\left(\int_{\overline{\Omega}_b}|\nabla v|^qd\mu+\int_{\Omega_b}|v|^qd\mu\right)-\int_{\Omega_a\cup\Omega_b} F(x,u,v)d\mu,
\end{eqnarray*}
where $\overline{\Omega}_a:=\Omega_a\cup\partial\Omega_a$ and $\overline{\Omega}_b:=\Omega_b\cup\partial\Omega_b$.
By a standard argument, we can verify that $J_\Omega\in C^1(W_\Omega,\R)$ and
\begin{eqnarray*}
\label{e7}
\langle J'_\Omega(u,v),(\phi_1,\phi_2)\rangle&=&\int_{\overline{\Omega}_a}|\nabla u|^{p-2} \Gamma(u,\phi_1)d\mu+\int_{\Omega_a}|u|^{p-2}u\phi_1d\mu
+\int_{\overline{\Omega}_b}|\nabla v|^{q-2} \Gamma(v,\phi_2)d\mu+\int_{\Omega_b}|v|^{q-2}v\phi_2d\mu\nonumber\\
&&-\int_{\Omega_a\cup\Omega_b} F_u(x,u,v)\phi_1d\mu-\int_{\Omega_a\cup\Omega_b} F_v(x,u,v)\phi_2d\mu
\end{eqnarray*}
for any $(u,v), (\phi_1,\phi_2)\in W_\Omega$. The corresponding Nehari manifold is
$$
\mathcal{N}_\Omega:=\{(u,v)\in W_\Omega\backslash\{(0,0)\}:\langle J'_\Omega(u,v),(u,v)\rangle=0\}.
$$
Let
\begin{eqnarray*}
\label{e8}
m_\Omega:=\inf_{(u,v)\in \mathcal{N}_\Omega}J_\Omega(u,v).
\end{eqnarray*}
\par
Similar to Theorem 1.1, the system (\ref{e5}) also has a ground state solution and the ground state solutions $\{(u_{\lambda}, v_{\lambda})\}$ of system (\ref{eq1})
converge to a ground state solution of system (\ref{e5}).

\vskip2mm
\noindent
{\bf Theorem 1.2.} {\it Let $G=(V, E)$ be a locally finite graph and it is connected, symmetric and there exists $\mu_{\min}>0$ such that $\mu(x)\geq \mu_{\min}$ for all $x\in \Omega$. For any $x\in V$, $\sum_{y\thicksim x}\omega_{xy} <C$, where $C$ is a positive constant. $\Omega_a, \Omega_b$ and $\Omega_a\cap\Omega_b$ are non-empty, connected and bounded domains in $V$. Assume $\frac{1}{\gamma}<\beta +\frac{1}{\beta} -1$, $(A_1)$, $(A_2)$ and $F:V\times\R\times\R\rightarrow\R$ satisfies $(F_1)-(F_4)$. Then system (\ref{e5}) has a ground state solution $(u_0, v_0)\in W_\Omega$. Besides for any sequence $\lambda_k\rightarrow+\infty$, up to a sequence, the ground state solutions family $\{(u_{\lambda_k}, v_{\lambda_k})\}$ of (\ref{eq1}) converge to a ground state solution of (\ref{e5}).}

\vskip2mm
\noindent
{\bf Remark 1.2.} Our works generalize those results in \cite{Han X L 2021} in some sense, and $(F_2)$ corresponding to the scalar case is weaker than $(f_2)$. Moreover, in \cite{Shao 2023}, Shao studied system (\ref{eq1}) with $p=q$ and  $p\geq2$ and the nonlinear term $F$ is independent on $x$ (see system (\ref{f1}) for details). However, our system (\ref{eq1}) allowed the fact that $p,q>1$ and $p,q$ are allowed to be unequal, and  the nonlinear term $F$ is allowed to depend on $x$.   It is easy to verify that the example in Remark 1.1 does not satisfy the condition $(F)$. Hence our result is different from that in \cite{Shao 2023}.

\vskip2mm
{\section{Preliminaries}}
  \setcounter{equation}{0}
  \par
In this section, we present some Sobolev embedding theorems on the locally finite graph and their proofs.
 \vskip2mm
\noindent
{\bf Lemma 2.1.} {\it Let $G=(V,E)$ be a locally finite graph and assume $p,q>1$. Then for any $\lambda \geq1$ and all $(\psi, \phi)\in W_\lambda$, there is
$$
\|\psi\|_{\infty, V}\leq d_1\|\psi\|_{W_\lambda(a)},\quad\|\phi\|_{\infty, V}\leq d_2\|\phi\|_{W_\lambda(b)},\quad\|(\psi,\phi)\|_{\infty, V}\leq d\|(\psi,\phi)\|_{W_\lambda},
$$
where $d_1=\left(\frac{1}{\mu_{\min} }\right)^\frac{1}{p}, d_2=\left(\frac{1}{\mu_{\min} }\right)^\frac{1}{q}$ and $d=\max\{d_1,d_2\}$.}
\vskip2mm
\noindent
{\bf Proof}. It follows from Lemma 2.6 in \cite{Han X L 2021} that  $\|\psi\|_{\infty, V}\leq d_1\|\psi\|_{W_\lambda(a)}$ and $\|\phi\|_{\infty, V}\leq d_2\|\phi\|_{W_\lambda(b)}$.  Furthermore,
$$
\|(\psi,\phi)\|_{\infty, V}=\max_{x\in V}|(\psi,\phi)|\leq\|\psi\|_{\infty, V}+\|\phi\|_{\infty, V} \leq d_1\|\psi\|_{W_\lambda(a)}+d_2\|\phi\|_{W_\lambda(b)}\leq \max\{d_1,d_2\}\|(\psi,\phi)\|_{W_\lambda}.
$$
The proof is complete.
\qed

 \vskip2mm
\noindent
{\bf Lemma 2.2.}  {\it Let $G=(V,E)$ be a locally finite graph and assume that $(A_1)$ and $(A_2)$ hold. Then for any $\lambda \geq1$, $W_\lambda(a)$ and $W_\lambda(b)$ are continuously embedded into $L^{\theta_1}(V)$ and $L^{\theta_2}(V)$ and  for any $\psi \in W_\lambda(a)$ and $\phi \in W_\lambda(b)$,
\begin{eqnarray}
\label{Eq11}
\|\psi\|_{\theta_1}\leq K_{(p,\theta_1)}\|\psi\|_{W_\lambda(a)},\;\;\;\|\phi\|_{\theta_2}\leq K_{(q,\theta_2)}\|\phi\|_{W_\lambda(b)},
\end{eqnarray}
where $\theta_1, \theta_2\in [1,+\infty)$,
\begin{eqnarray*}
&&K_{(p,\theta_1)}=
 \begin{cases}
d_1^{\frac{\theta_1-p}{\theta_1}},\;\;\;\;\;\;\;\;\;\;\;\;\;\;\;\;\;\;\;\;\;\;\;\;\;\;\;\;\;\;\;\;\;\;\;\;\;\;\;\;\;\;\;\;\;\;\;\;\;\;\;\;\;\;\;\;\;\;\;\;\;\;\;\;
\;\;\;\;\;\;\;\;\;\;\;\;\;\;\;\;\;\mbox{ $\theta_1\geq p$},\\
\min\left\{d_1^{\frac{\theta_1-1}{\theta_1}}\|(a+1)^{-1}\|^{\frac{1}{p\theta_1}}_{\frac{1}{p-1}}, d_1^\frac{p(\theta_1-1)}{\theta_1}\|(a+1)^{-1}\|^{\frac{1}{p}}_{\frac{1}{p-1}}\right\},\;\;\;\;\;\;\;\; \mbox{$1\leq\theta_1< p$}, \\
\end{cases}
\end{eqnarray*}
\begin{eqnarray*}
&&K_{(q,\theta_2)}=
 \begin{cases}
d_2^{\frac{\theta_2-q}{\theta_2}},\;\;\;\;\;\;\;\;\;\;\;\;\;\;\;\;\;\;\;\;\;\;\;\;\;\;\;\;\;\;\;\;\;\;\;\;\;\;\;\;\;\;\;\;\;\;\;\;\;\;\;\;\;\;\;\;\;\;\;\;\;\;\;\;
\;\;\;\;\;\;\;\;\;\;\;\;\;\;\;\;\;\mbox{ $\theta_2\geq q$},\\
\min\left\{d_2^{\frac{\theta_2-1}{\theta_2}}\|(b+1)^{-1}\|^{\frac{1}{q\theta_2}}_{\frac{1}{q-1}}, d_2^\frac{q(\theta_2-1)}{\theta_2}\|(b+1)^{-1}\|^{\frac{1}{q}}_{\frac{1}{q-1}}\right\},\;\;\;\;\;\;\;\;\mbox{$1\leq\theta_2< q$}.\\
\end{cases}
\end{eqnarray*}
Moreover, $W_\lambda$ is continuously embedded into $L(V):=L^{\theta_1}(V)\times L^{\theta_2}(V)$ for all $\theta_1, \theta_2\in [1,+\infty)$ and
\begin{eqnarray}
\label{Z3}
\|(\psi,\phi)\|_{L(V)}\leq \max\{K_{(p,\theta_1)},K_{(q,\theta_2)}\}\|(\psi,\phi)\|_{W_\lambda},
\end{eqnarray}

}
\vskip2mm
\noindent
{\bf Proof.}\ \ The proof is similar to Lemma 2.6 in \cite{Han X L 2021} with some slight modifications. According to Lemma 2.1, when $\theta_1\geq p$, we have
$$
\|\psi\|_{\theta_1}^{\theta_1}
  =\int_V|\psi|^{\theta_1-p}\cdot |\psi|^{p}d\mu
\leq\|\psi\|_{\infty, V}^{\theta_1-p}\int_V (|\nabla  \psi|^p+(\lambda a+1)|\psi(x)|^p)d\mu
\leq d_1^{\theta_1-p}\|\psi\|_{W_\lambda(a)}^{\theta_1-p}\cdot\|\psi\|_{W_\lambda(a)}^{p}
  = d_1^{\theta_1-p}\|\psi\|_{W_\lambda(a)}^{\theta_1 }.
$$
By $(A_2)$ and noticing that $\lambda\geq1$, we have
\begin{eqnarray}\label{a2}
\int_V|\psi|d\mu&=&\int_V(\lambda a+1)^{-\frac{1}{p}}(\lambda a+1)^{\frac{1}{p}}|\psi|d\mu\nonumber\\
&   \leq&\left(\int_V(\lambda a+1)^{-\frac{1}{p}\times\frac{p}{ p-1}}d\mu\right)^{\frac{p-1}{p}} \left(\int_V(\lambda a+1)^{\frac{1}{p}\times p}|\psi|^{p}d\mu\right)^{\frac{1}{p}}\nonumber\\
&\leq&\|(a+1)^{-1}\|_{\frac{1}{ p-1}}^{\frac{1}{p}}\|\psi\|_{W_\lambda(a)}.
\end{eqnarray}
On one hand, when $1\leq\theta_1< p$, from (\ref{a2}), we have
\begin{eqnarray}
\label{a1}
\|\psi\|_{\theta_1}^{\theta_1}
 =\int_V|\psi|^{\theta_1-1}\cdot |\psi|d\mu
\leq\|\psi\|_{\infty, V}^{\theta_1-1}\int_V |\psi| d\mu
\leq d_1^{\theta_1-1}\|(a+1)^{-1}\|_{\frac{1}{ p-1}}^{\frac{1}{p}}\|\psi\|_{W_\lambda(a)}^{\theta_1 }.
\end{eqnarray}
On the other hand, by Lemma 2.6 in \cite{Han X L 2021}, $\|\psi\|_{\theta_1} \leq d_1^\frac{p(\theta_1-1)}{\theta_1}\|(a+1)^{-1}\|^{\frac{1}{p}}_{\frac{1}{p-1}}\|\psi\|_{W_\lambda(a)}$. So combining (\ref{a1}), there is
$$
\|\psi\|_{\theta_1} \leq  \min\left\{d_1^{\frac{\theta_1-1}{\theta_1}}\|(a+1)^{-1}\|^{\frac{1}{p\theta_1}}_{\frac{1}{p-1}}, d_1^\frac{p(\theta_1-1)}{\theta_1}\|(a+1)^{-1}\|^{\frac{1}{p}}_{\frac{1}{p-1}}\right\}\|\psi\|_{W_\lambda(a)}.
$$
Similarly, it is easy to prove that the second inequality also holds in (\ref{Eq11}).Thus,
$$
\|(\psi,\phi)\|_{L(V)}=\|\psi\|_{\theta_1}+\|\phi\|_{\theta_2}\leq K_{(p,\theta_1)}\|\psi\|_{W_\lambda(a)}+K_{(q,\theta_2)}\|\phi\|_{W_\lambda(b)} \leq \max\{K_{(p,\theta_1)},K_{(q,\theta_2)}\}\|(\psi,\phi)\|_{W_\lambda}.
$$
The proof is complete.
\qed

 \vskip2mm
\noindent
{\bf Lemma 2.3.} {\it Let $G=(V, E)$ be a locally finite graph and assume that $(A_1)$ and $(A_2)$ hold. Then for any given $\lambda \geq1$ and any bounded sequence $(\psi_k, \phi_k) \in W_\lambda$, there exists $(\psi,\phi)\in W_\lambda$ such that, up to subsequence,
\begin{eqnarray*}
 \begin{cases}
  (\psi_k,\phi_k)\rightharpoonup (\psi,\phi),\;\;\;\;\;\;\;\;\;\;\;\;\;\;\;\;\;\;\;\;\;\mbox{in}\; W_\lambda,\;\;\mbox{as $k \rightarrow \infty$,}\\
  \psi_k(x)\rightarrow \psi(x),\phi_k(x)\rightarrow\phi(x),\;\;\;\;\; \forall x\in V,\; \mbox{as $k\rightarrow \infty$},\\
  \psi_k \rightarrow \psi,\;\;\;\;\;\;\;\;\;\;\;\;\;\;\;\;\;\;\;\;\;\;\;\;\;\;\;\;\;\;\;\;\;\;\;\;\mbox{in}\; L^{\theta_1}(V), \;\;\forall \theta_1\in[1,+\infty],\; \mbox{as $k\rightarrow \infty$},\\
  \phi_k \rightarrow \phi,\;\;\;\;\;\;\;\;\;\;\;\;\;\;\;\;\;\;\;\;\;\;\;\;\;\;\;\;\;\;\;\;\;\;\;\;\;\mbox{in}\; L^{\theta_2}(V), \;\;\forall \theta_2\in[1,+\infty],\; \mbox{as $k\rightarrow \infty$}.
   \end{cases}
\end{eqnarray*}
}
\vskip0mm
\noindent
{\bf Proof.}\ \ Since $W_\lambda$ is a reflexive Banach space. Thus for any bounded sequence ${(\psi_k,\phi_k)}\in W_\lambda$, we get that, up to a subsequence, $(\psi_k,\phi_k)\rightharpoonup (\psi,\phi)$ in $W_\lambda$.
In addition, by Lemma 2.2, we can get $\psi_k\rightharpoonup \psi$ in $L^{\theta_1}(V)$ and $\phi_k\rightharpoonup\phi$ in $L^{\theta_2}(V)$. The remaining proof is similar to Lemma 2.6 in \cite{Han X L 2021} with substituting $\psi_k$ and $\phi_k$ for $u_k$ in \cite{Han X L 2021}, respectively. We omit it here.\qed

 \vskip2mm
\noindent
{\bf Lemma 2.4.} {\it $W_{\Omega_a}$ and $W_{\Omega_b}$ are compactly embedded into $L^{\theta_1}(\Omega_a)$ and $L^{\theta_2}(\Omega_b)$ for any $\theta_1, \theta_2\in [1,+\infty]$, respectively,  and  for any $\psi \in W_\Omega$,
\begin{eqnarray}
\label{e14}
\|\psi\|_{\theta_1,\Omega_a}\leq K_{(p,\theta_1)}^*\|\psi\|_{W_{\Omega_a}},\;\;\|\phi\|_{\theta_2,\Omega_b}\leq K_{(q,\theta_2)}^*\|\phi\|_{W_{\Omega_b}},\;\;\|(\psi,\phi)\|_{\theta, \Omega}\leq K^*\|(\psi,\phi)\|_{W_{\Omega}},
\end{eqnarray}
where
\begin{eqnarray*}
K_{(p,\theta_1)}^*=\left(\sum_{x\in\Omega_a}\mu(x)\right)^{\frac{1}{\theta_1}}\left(\frac{1}{\mu_{\min}}\right)^{\frac{1}{p}},
K_{(q,\theta_2)}^*=\left(\sum_{x\in\Omega_b}\mu(x)\right)^{\frac{1}{\theta_2}}\left(\frac{1}{\mu_{\min}}\right)^{\frac{1}{q}},
 K^*=\max\left\{K_{(p,\theta_1)}^*, K_{(q,\theta_2)}^*\right\}.
\end{eqnarray*}
Moreover, $W_{\Omega_a}(W_{\Omega_b})$ is pre-compact, namely, if $\psi_k(\phi_k)$ is bounded in $W_{\Omega_a}(W_{\Omega_b})$, then up to a subsequence, there exists some $\psi \in W_{\Omega_a}(\phi \in W_{\Omega_b})$ such that $\psi_k \rightarrow \psi$ in $W_{\Omega_a}(\phi_k \rightarrow \phi$ in $W_{\Omega_b})$.}

\vskip2mm
\noindent
{\bf Proof.}\ \  The proof of (\ref{e14}) is similar to Lemma 2.1 and Lemma 2.2. In fact,
\begin{eqnarray}
\label{t2}
         \|\psi\|_{W_{\Omega_a}}^p
 =   \int_{\overline{\Omega}_a}|\nabla  \psi|^pd\mu+\int_{{\Omega}_a}|\psi(x)|^pd\mu
 \geq   \sum_{x\in  \Omega_a }\mu(x) |\psi(x)|^p
 \geq   \mu_{\min}\sum_{x\in \Omega_a}|\psi(x)|^p
 \geq   \mu_{\min}\|\psi\|^p_{\infty, \Omega_a}.
\end{eqnarray}
From (\ref{t2}), we have
\begin{eqnarray*}
\|\psi\|_{\theta_1,\Omega_a}^{\theta_1}
 =\sum_{x\in\Omega_a}\mu(x)|\psi|^{\theta_1}
\leq\sum_{x\in\Omega_a}\mu(x) \|\psi\|_{\infty, \Omega_a}^{\theta_1}
\leq\sum_{x\in\Omega_a}\mu(x)\left(\frac{1}{\mu_{\min}}\right)^{\frac{\theta_1}{p}}\|\psi\|_{W_{\Omega_a}}^{\theta_1}.
\end{eqnarray*}
The proof of the first inequality in (\ref{e14}) is completed. Similarly, we can also easily complete the proofs of the other two inequalities in (\ref{e14}). Since $\Omega_a(\Omega_b)$ is a finite set in $V$, $W_{\Omega_a}(W_{\Omega_b})$ and $L^{\theta_1}(\Omega_a)(L^{\theta_2}(\Omega_b))$ are finite dimensional spaces. Hence, $W_{\Omega_a}(W_{\Omega_b})$ is pre-compact. Then for any bounded sequence $\{(\psi_k,\phi_k)\}\in W_\Omega$, there is  subsequence, still denoted by $\{(\psi_k,\phi_k)\}\in W_\Omega$, such that  $\|\psi_k-\psi\|_{W_{\Omega_a}}\rightarrow 0$ and $\|\phi_k-\phi\|_{W_{\Omega_b}}\rightarrow 0$. Furthermore,  by (\ref{e14}),  we can get $\|\psi_k-\psi\|_{\theta_1, \Omega_a}\rightarrow 0$ and $\|\phi_k-\phi\|_{\theta_2,\Omega_b}\rightarrow 0$. So $W_{\Omega_a}(W_{\Omega_b})$ is compactly embedded into $L^{\theta_1}(\Omega_a)(L^{\theta_2}(\Omega_b))$ for any $\theta_1, \theta_2\in [1,+\infty]$.
\qed

\vskip2mm
{\section{The existence of ground state solutions}}
  \setcounter{equation}{0}
  \par
In this section we prove Theorem 1.1 by the Nehari manifold method and the Lagrange multiplier rule.
 \vskip2mm
\noindent
{\bf Definition 3.1.} Suppose that $(u,v)\in W_\lambda$. If for any $(\phi_1, \phi_2)\in W_\lambda$, there holds
\begin{eqnarray*}
&&\int_V(|\nabla u|^{p-2}\Gamma(u,\phi_1)+(\lambda a+1)|u|^{p-2}u\phi_1)d\mu
+\int_V(|\nabla v|^{q-2}\Gamma(v,\phi_2)+(\lambda b+1)|v|^{q-2}v\phi_2)d\mu\nonumber\\
&=&\int_V F_u(x,u,v)\phi_1d\mu+\int_V F_v(x,u,v)\phi_2d\mu,
\end{eqnarray*}
then $(u,v)$ is called a weak solution of (\ref{eq1}).

  \vskip2mm
\noindent
{\bf Definition 3.2.} Suppose that $(u,v)\in W_\Omega$. If for any $(\phi_1, \phi_2)\in W_\Omega$, there holds
\begin{eqnarray*}
&&\int_{\overline{\Omega}_a}|\nabla u|^{p-2}\Gamma(u,\phi_1)d\mu+\int_{\Omega_a}|u|^{p-2}u\phi_1d\mu
+\int_{\overline{\Omega}_b}|\nabla v|^{q-2}\Gamma(v,\phi_2)d\mu+\int_{\Omega_b}|v|^{q-2}v\phi_2d\mu\nonumber\\
&=&\int_{\Omega_a\cup\Omega_b} F_u(x,u,v)\phi_1d\mu+\int_{\Omega_a\cup\Omega_b} F_v(x,u,v)\phi_2d\mu,
\end{eqnarray*}
then $(u,v)$ is called a weak solution of (\ref{e5}).

  \vskip2mm
\noindent
{\bf Lemma 3.1.} {\it If $(u,v)\in W_\lambda$ is a weak solution of (\ref{eq1}), $(u,v)$ is also a point-wise solution of (\ref{eq1}).}
\vskip2mm
\noindent
{\bf Proof.}\ \ The proof is standard (for example, see \cite{Yang P 2023}). For completeness, we also present it here. Since $(u,v)\in W_\lambda$ is a weak solution of (\ref{eq1}), for any $(\phi,\psi)\in W_\lambda$, there holds
\begin{eqnarray*}
&&\int_V(|\nabla u|^{p-2}\Gamma(u ,\phi)+(\lambda a+1)|u|^{p-2}u\phi)d\mu
+\int_V(|\nabla v|^{q-2}\Gamma(v ,\psi)+(\lambda b+1)|v|^{q-2}v\psi)d\mu\nonumber\\
&=&\int_V F_u(x,u,v)\phi d\mu+\int_V F_v(x,u,v)\psi d\mu.
\end{eqnarray*}
From Lemma 2.1 in \cite{Han X L 2021}, we have
\begin{eqnarray}
\label{f7}
\int_V |\nabla u|^{p-2}\Gamma(u ,\phi) d\mu=\int_V(-\Delta_p u)\phi d\mu,\;\;\;\;\;\;\;\int_V |\nabla v|^{q-2}\Gamma(v ,\psi) d\mu=\int_V(-\Delta_q v)\psi d\mu.
\end{eqnarray}
Then we get
\begin{eqnarray}\label{EQ2}
&&\int_V((-\Delta_p u)\phi+ (\lambda a+1)|u|^{p-2}u \phi) d\mu+\int_V((-\Delta_q v)\psi+ (\lambda b+1)|v|^{q-2}v \psi )d\mu\nonumber\\
&=&\int_V F_u(x,u,v) \phi d\mu+\int_V F_v(x,u,v) \psi d\mu.
\end{eqnarray}
If for any $x_0 \in V$, we take the test function $(\phi(x),\psi(x)): V\rightarrow \R$ in (\ref{EQ2}) with

$$
(\phi(x), \psi(x))=
 \begin{cases}
  (1,0)\;\;\;\; x=x_0,\\
  (0,0)\;\;\;\; x\neq x_0,
   \end{cases}
,\;\;\;
(\phi(x), \psi(x))=
 \begin{cases}
  (0,1)\;\;\;\; x=x_0,\\
  (0,0)\;\;\;\; x\neq x_0,
   \end{cases}
$$
respectively, then
\begin{eqnarray*}
 \begin{cases}
  -\Delta_p u(x_0)+(\lambda a+1)|u(x_0)|^{p-2}u(x_0)=F_u(x,u(x_0),v(x_0)),\;\;\;\;\hfill x_0\in V,\\
  -\Delta_q v(x_0)+(\lambda b+1)|v(x_0)|^{q-2}v(x_0)=F_v(x,u(x_0),v(x_0)),\;\;\;\;\hfill x_0\in V.\\
   \end{cases}
\end{eqnarray*}
Since $x_0$ is arbitrary, we conclude that $(u,v)$ is a point-wise solution of (\ref{eq1}).
\qed
\vskip2mm
\noindent
{\bf Lemma 3.2.} {\it If $(u,v) \in W_\Omega$ is a weak solution of (\ref{e5}), then $(u,v)$ is also a point-wise solution of (\ref{e5}).}
\vskip2mm
\noindent
{\bf Proof.}\ \ The proof is similar to Lemma 3.1. Since $(u,v)\in W_\Omega$ is a weak solution of (\ref{e5}), for any $(\phi,\psi)\in C_c(\Omega_a)\times C_c(\Omega_b)$, there holds
\begin{eqnarray*}
&&\int_{\overline{\Omega}_a}|\nabla u|^{p-2}\Gamma(u,\phi)d\mu+\int_{\Omega_a}|u|^{p-2}u\phi d\mu+\int_{\overline{\Omega}_b}|\nabla v|^{q-2}\Gamma(v,\psi)d\mu+\int_{\Omega_b}|v|^{q-2}v\psi d\mu\nonumber\\
&=&\int_{\Omega_a\cup\Omega_b} F_u(x,u,v)\phi d\mu+\int_{\Omega_a\cup\Omega_b} F_v(x,u,v)\psi d\mu.
\end{eqnarray*}
Similar to (\ref{f7}), we also have
\begin{eqnarray*}
\int_{\overline{\Omega}_a} |\nabla u|^{p-2}\Gamma(u ,\phi) d\mu=\int_{\overline{\Omega}_a}(-\Delta_p u)\phi d\mu,\;\;\;\;\;\;\;\int_{\overline{\Omega}_b} |\nabla v|^{q-2}\Gamma(v ,\psi) d\mu=\int_{\overline{\Omega}_b}(-\Delta_q v)\psi d\mu.
\end{eqnarray*}
and then
\begin{eqnarray}\label{Z8}
&&\int_{\overline{\Omega}_a}(-\Delta_p u)\phi d\mu+\int_{\Omega_a}|u|^{p-2}u\phi d\mu+\int_{\overline{\Omega}_b}(-\Delta_q v)\psi d\mu+\int_{\Omega_b}|v|^{q-2}v\psi d\mu\nonumber\\
&=&\int_{\Omega_a\cup\Omega_b} F_u(x,u,v)\phi d\mu+\int_{\Omega_a\cup\Omega_b} F_v(x,u,v)\psi d\mu.
\end{eqnarray}
If for any given $x_1 \in \Omega_a$, we take the test function $(\phi(x),\psi(x)): \Omega_a\times \Omega_b \rightarrow \R\times \R$ in (\ref{Z8}) with
\begin{eqnarray*}
(\phi(x), \psi(x))=
 \begin{cases}
  (1,0)\;\;\;\; x=x_1,\\
  (0,0)\;\;\;\; x\neq x_1,
   \end{cases}
\end{eqnarray*}
and for any given $x_2 \in \Omega_b$, we take the test function $(\phi(x),\psi(x)): \Omega_a\times \Omega_b \rightarrow \R\times \R$ in (\ref{Z8}) with
\begin{eqnarray*}
(\phi(x), \psi(x))=
 \begin{cases}
  (0,1)\;\;\;\; x=x_2,\\
  (0,0)\;\;\;\; x\neq x_2,
   \end{cases}
\end{eqnarray*}
respectively, then
\begin{eqnarray*}
 \begin{cases}
  -\Delta_p u(x_1)+|u(x_1)|^{p-2}u(x_1)=F_u(x,u(x_1),v(x_1)),\\
  -\Delta_q v(x_2)+|v(x_2)|^{q-2}v(x_2)=F_v(x,u(x_2),v(x_2)).\\
   \end{cases}
\end{eqnarray*}
Since $W_{\Omega_a}(W_{\Omega_b})$ is the completion of $C_c(\Omega_a)(C_c(\Omega_b))$,  $u(x)=0$ on $\partial \Omega_a$ and $v(x)=0$ on $\partial \Omega_b$. Finally, by the arbitrary of $x_1$ and $x_2$,  we complete the proof.
\qed

 \vskip2mm
\noindent
{\bf Lemma 3.3.} {\it If $(F_1)$ and $(F_4)$ hold and $(u,v) \in W_\lambda\backslash\{(0,0)\}$, then there exists a unique $t_0>0$ such that $(t_0u, t_0v)\in \mathcal{N_\lambda}$ and then $ \mathcal{N_\lambda}$ is non-empty.}
\vskip2mm
\noindent
{\bf Proof.}\ \ We define a fibering map $g(t):=J_\lambda(tu,tv)$ on $\R$ for all $t>0$. Clearly, we have
\begin{eqnarray*}
g'(t)=0 \Longleftrightarrow \langle J'_\lambda(tu,tv),(u,v)\rangle=0
\Longleftrightarrow \langle J'_\lambda(tu,tv),(tu,tv)\rangle=0
\Longleftrightarrow (tu,tv)\in \mathcal{N_\lambda}.
\end{eqnarray*}
Since $t>0$ and $\beta=\max\{p,q\}$, then
\begin{eqnarray}
\label{t1}
&&\frac{g'(t)}{t^{\beta-1}}=0
 \Longleftrightarrow g'(t)=0
  \Longleftrightarrow(tu,tv)\in \mathcal{N_\lambda}\nonumber\\
  &\Longleftrightarrow& t^{p-\beta}\|u\|^p_{W_\lambda(a)}+t^{q-\beta}\|v\|^q_{W_\lambda(b)}=t^{-\beta}\int_V[F_{tu}(x,tu,tv)tu+F_{tv}(x,tu,tv)tv]d\mu.
\end{eqnarray}
Without loss of generality, we let $\beta=p$. Then we infer that
\begin{eqnarray}
\label{b2}
\|u\|^p_{W_\lambda(a)}=t^{1-p}\int_V[F_{tu}(x,tu,tv)u+F_{tv}(x,tu,tv)v]d\mu-t^{q-p}\|v\|^q_{W_\lambda(b)}.
\end{eqnarray}
Let $H_1(t)=t^{1-p}\int_V[F_{tu}(x,tu,tv)u+F_{tv}(x,tu,tv)v]d\mu-t^{q-p}\|v\|^q_{W_\lambda(b)}$. Then from $(F_4)$, we have
\begin{eqnarray*}
&&H_1'(t)\\
&=&(1-p)t^{-p}\int_V\left[F_{tu}(x,tu,tv)u+F_{tv}(x,tu,tv)v\right]d\mu\\
&&+t^{1-p}\int_V\left[F_{tu,tu}(x,tu,tv)u^2+F_{tu,tv}(x,tu,tv)uv+F_{tv,tv}(x,tu,tv)v^2+F_{tv,tu}(x,tu,tv)v u\right]d\mu\\
&&-(q-p)t^{q-p-1}\|v\|^q_{W_\lambda(b)}\\
&=&(1-p)t^{-p}\int_V\left[F_{tu}(x,tu,tv)u+F_{tv}(x,tu,tv)v\right]d\mu\\
&&+t^{-p}\int_V\left[F_{tu,tu}(x,tu,tv)tu^2+F_{tv,tv}(x,tu,tv)tv^2\right]d\mu\\
&&+t^{1-p}\int_V\left[F_{tu,tv}(x,tu,tv)uv+F_{tv,tu}(x,tu,tv)v u\right]d\mu -(q-p)t^{q-p-1}\|v\|^q_{W_\lambda(b)}\\
&>&(1-p)t^{-p}\int_V\left[F_{tu}(x,tu,tv)u+F_{tv}(x,tu,tv)v\right]d\mu\\
&&+(p-1)t^{-p}\int_V\left[F_{tu}(x,tu,tv)u+F_{tv}(x,tu,tv)v\right]d\mu\\
&&+t^{1-p}\int_V\left[F_{tu,tv}(x,tu,tv)uv+F_{tv,tu}(x,tu,tv)v u\right]d\mu-(q-p)t^{q-p-1}\|v\|^q_{W_\lambda(b)}\\
&\geq&0.
\end{eqnarray*}
Thus $H_1(t)$ is strictly increasing in $(0, +\infty)$. Hence by observing (\ref{b2}), there exists a unique $t_0>0$ such that (\ref{t1}) holds and $(t_0u,t_0v)\in\mathcal{N_\lambda}$. Therefore $\mathcal{N_\lambda}$ is non-empty.
\qed

 \vskip2mm
\noindent
{\bf Lemma 3.4.} {\it If $(F_1)$ and $(F_4)$ hold, then $J_\lambda(tu,tv)\leq J_\lambda(u,v)$ for any given ${(u,v) \in\mathcal{N_\lambda}}$ and all $t>0$.}
\vskip2mm
\noindent
{\bf Proof.}\ \ Since ${(u,v) \in\mathcal{N_\lambda}}$, we have $g'(1)=0$ where $g(\cdot)$ is introduced in the proof of Lemma 3.3, and 1 is the unique critical point of $g(\cdot).$ Then
\begin{eqnarray*}
g(t)=J_\lambda(tu,tv)=\frac{t^p}{p}\| u\|^p_{W_\lambda(a)}+\frac{t^q}{q}\|v\|^q_{W_\lambda(b)}-\int_V F(x,tu,tv)d\mu.
\end{eqnarray*}
Then by $(F_4)$, we get
\begin{eqnarray*}
g'(t) &=  &t^{p-1}\| u\|^p_{W_\lambda(a)}+ t^{q-1}\|v \|^q_{W_\lambda(b)}-\int_V [F_{tu}(x,tu,tv)u+F_{tv}(x,tu,tv)v]d\mu.\\
g''(t)&=  &(p-1)t^{p-2}\| u\|^p_{W_\lambda(a)}+(q-1)t^{q-2}\|v \|^q_{W_\lambda(b)}- \int_V [F_{tu,tu}(x,tu,tv)u^2\\
&&+F_{tv,tv}(x,tu,tv)v^2+F_{tu,tv}(x,tu,tv)uv+F_{tv,tu}(x,tu,tv)vu]d\mu\\
&  <  & (p-1)t^{p-2}\| u\|^p_{W_\lambda(a)}+(q-1)t^{q-2}\|v \|^q_{W_\lambda(b)} - \int_V [(\beta-1)\frac{1}{t}F_{tu}(x,tu,tv)u\\
&&+(\beta-1)\frac{1}{t}F_{tv}(x,tu,tv)v+F_{tu,tv}(x,tu,tv)uv+F_{tv,tu}(x,tu,tv)vu]d\mu.
\end{eqnarray*}
Let $t=1$. By $(F_4)$, we have
\begin{eqnarray}
\label{Z6}
g''(1)& <  & (p-1) \| u\|^p_{W_\lambda(a)}+(q-1) \|v \|^q_{W_\lambda(b)} - \int_V [(\beta-1) F_u(x, u, v)u\nonumber\\
&&+(\beta-1) F_v(x, u, v)v+F_{uv}(x, u, v)uv+F_{vu}(x, u, v)vu]d\mu\nonumber \\
&  \leq  &  (\beta-1) \left(\| u\|^p_{W_\lambda(a)}+ \|v \|^q_{W_\lambda(b)} \right)- \int_V [(\beta-1) F_u(x, u, v)u\nonumber\\
&&+(\beta-1) F_v(x, u, v)v+F_{uv}(x, u, v)uv+F_{vu}(x, u, v)vu]d\mu\nonumber \\
&  =  & (\beta-1)\int_V [F_u(x, u, v)u+F_v(x, u, v)v]d\mu- \int_V [(\beta-1) F_u(x, u, v)u\nonumber\\
&&+(\beta-1) F_v(x, u, v)v+F_{uv}(x, u, v)uv+F_{vu}(x, u, v)vu]d\mu \nonumber\\
&  \leq &  0.
\end{eqnarray}
\par
From (\ref{Z6}) and the fact that 1 is the unique critical point of $g(t)$, we obtain that $t=1$ is the maximizer of $g(\cdot)$. Hence $J_\lambda(tu,tv)\leq J_\lambda(u,v)$ for all $t>0$.
\qed

\vskip2mm
\noindent
{\bf Lemma 3.5.} {\it If $(F_3)$ holds, there exists a constant $\eta >0$ such that $m_\lambda=\inf\limits_{(u,v) \in\mathcal{N_\lambda}}J_\lambda(u,v)\geq \eta >0$ .}
\vskip2mm
\noindent
{\bf Proof.}\ \ By $(F_3)$ and Appendix A.2, we can get
\begin{eqnarray}
\label{t4}
|F(x,s,t)|
\leq \frac{2C_1(x)}{p}|s|^p+\left(\frac{C_1(x)(p-1)}{p}+\frac{C_1(x)}{q}\right)|t|^q
+\frac{2C_2(x)}{r_1}|s|^{r_1}+\left(\frac{C_2(x)(r_1-1)}{r_1}+\frac{C_2(x)}{r_2}\right)|t|^{r_2}
\end{eqnarray}
for all $x\in V$. For $(u,v)\in W_\lambda\backslash\{(0,0)\}$ with $\|(u,v)\|_{W_\lambda}<1$, from (\ref{t4}) and Lemma 2.2, we have
\begin{eqnarray*}
\label{Z5}
&&J_\lambda( u, v)\nonumber\\
&\geq&\frac{1}{p}\| u\|^p_{W_\lambda(a)}+\frac{1}{q}\|v\|^q_{W_\lambda(b)}-\int_V\frac{2C_1(x)}{p}| u|^pd\mu-\int_V \left(\frac{C_1(x)(p-1)}{p}+\frac{C_1(x)}{q}\right) |v|^q d\mu\nonumber\\
&&-\int_V\frac{2C_2(x)}{r_1} |u | ^{r_1} d\mu-\int_V \left(\frac{C_2(x)(r_1-1)}{r_1}+\frac{C_2(x)}{r_2}\right) |v | ^{r_2}d\mu\nonumber\\
&\geq&\frac{1}{p}\| u\|^p_{W_\lambda(a)}+\frac{1}{q}\|v\|^q_{W_\lambda(b)}-\frac{2\|C_1\|_{\infty, V}}{p}\| u\|_p^p- \left(\frac{\|C_1\|_{\infty, V}(p-1)}{p}+\frac{\|C_1\|_{\infty, V}}{q}\right) \|v\|_q^q \nonumber\\
&&-\frac{2\|C_2\|_{\infty, V}}{r_1} \|u \|_ {r_1}^{r_1} - \left(\frac{\|C_2\|_{\infty, V}(r_1-1)}{r_1}+\frac{\|C_2\|_{\infty, V}}{r_2}\right) \|v \|_{r_2} ^{r_2}\nonumber\\
&\geq&\left(\frac{1}{p}-\frac{2\|C_1\|_{\infty, V}}{p}-\frac{2\|C_2\|_{\infty, V}}{r_1}K_{(p,r_1)}^{r_1}\|u\|_{W_\lambda(a)}^{r_1-p}\right)\|u\|^p_{W_\lambda(a)}\nonumber\\
&&+\left[\frac{1}{q}-\left(\frac{\|C_1\|_{\infty, V}(p-1)}{p}+\frac{\|C_1\|_{\infty,V}}{q}\right)-\left(\frac{\|C_2\|_{\infty,V}(r_1-1)}{r_1}+\frac{\|C_2\|_{\infty,V}}{r_2}\right)K_{(q,r_2)}^{r_2}\|v\|_{W_\lambda(b)}^{r_2-q}\right]\|v\|^q_{W_\lambda(b)}\\
&\geq&\frac{1}{2^{\beta-1}}\|(u,v)\|^\beta_{W_\lambda}\min\left\{\frac{1}{p}-\frac{2\|C_1\|_{\infty, V}}{p}-\frac{2\|C_2\|_{\infty, V}}{r_1}K_{(p,r_1)}^{r_1}\|(u,v)\|_{W_\lambda}^{r_1-p},\right.\nonumber\\
&&\left.\frac{1}{q}-\left(\frac{\|C_1\|_{\infty, V}(p-1)}{p}+\frac{\|C_1\|_{\infty,V}}{q}\right)
-\left(\frac{\|C_2\|_{\infty,V}(r_1-1)}{r_1}+\frac{\|C_2\|_{\infty,V}}{r_2}\right)K_{(q,r_2)}^{r_2}\|(u,v)\|_{W_\lambda}^{r_2-q}\right\}
\end{eqnarray*}
Let
$$\rho \in \left(0,\min\left\{1,\left(\frac{\frac{1}{p}-\frac{2\|C_1\|_{\infty,V}}{p}}{\frac{2\|C_2\|_{\infty,V}}{r_1}K_{(p,r_1)}^{r_1}}\right)^{\frac{1}{r_1-p}}, \left(\frac{\frac{1}{q}-\frac{\|C_1\|_{\infty,V}(p-1)}{p}-\frac{\|C_1\|_{\infty,V}}{q}}{(\frac{\|C_2\|_{\infty,V}(r_1-1)}{r_1}+\frac{\|C_2\|_{\infty,V}}{r_2})K_{(q,r_2)}^{r_2}}\right)^{\frac{1}{r_2-q}}\right\}\right)
$$
where $p<r_1$, $q<r_2$ and $C_1(x)\leq \frac{1}{1+\beta}<\frac{1}{2}$, and
\begin{eqnarray*}
&&\eta\\
&=&\frac{1}{2^{\beta-1}}\rho^\beta\min\left\{\frac{1}{p}-\frac{2\|C_1\|_{\infty, V}}{p}-\frac{2\|C_2\|_{\infty, V}}{r_1}K_{(p,r_1)}^{r_1}\rho^{r_1-p},\right.\nonumber\\
&&\left.\frac{1}{q}-\left(\frac{\|C_1\|_{\infty, V}(p-1)}{p}+\frac{\|C_1\|_{\infty,V}}{q}\right)
-\left(\frac{\|C_2\|_{\infty,V}(r_1-1)}{r_1}+\frac{\|C_2\|_{\infty,V}}{r_2}\right)K_{(q,r_2)}^{r_2}\rho^{r_2-q}\right\}.
\end{eqnarray*}

Then
\begin{eqnarray*}
J_\lambda( u, v) \geq \eta>0\;\;\;\mbox{for\;all}\;\|(u,v)\|_{W_\lambda}=\rho.
\end{eqnarray*}
For any given $(u,v) \in\mathcal{N_\lambda}\subset W_\lambda$, there exists $\tau_{\lambda,u,v}>0$ such that $\tau_{\lambda,u,v}\|(u,v)\|_{W_\lambda}=\rho$. Then by Lemma 3.4, we have
\begin{eqnarray*}
\label{t3}
&&J_\lambda( u, v) \geq J_\lambda(\tau_{\lambda,u,v} u,\tau_{\lambda,u,v} v) \geq \eta>0\\
&\Rightarrow& m_\lambda=\inf\limits_{(u,v) \in\mathcal{N_\lambda}}J_\lambda(u,v)>0.
\end{eqnarray*}
Thus we finish the proof.
\qed

 \vskip2mm
\noindent
{\bf Lemma 3.6.} {\it If $(F_2)$ holds, then $J_\lambda|_\mathcal{N_\lambda}$ is coercive.}
\vskip2mm
\noindent
{\bf Proof.}\ \ We prove the inverse negative proposition of coercive. That is, if $\{(u_k,v_k)\}\subseteq\mathcal{N_\lambda}$ and $J_\lambda( u_k, v_k)\leq M$ for some $M>0$ and all $k \in \mathbb{N}$, then $\{(u_k,v_k)\}\subseteq \mathcal{N_\lambda}$ is bounded.
\par
Indeed, from $(F_2)$, we have
\begin{eqnarray*}
\label{t11}
M
&\geq&\frac{1}{p}\| u_k\|^p_{W_\lambda(a)}+\frac{1}{q}\| v_k\|^q_{W_\lambda(b)}-\int_V  F(x,u_k,v_k)d\mu\nonumber\\
&\geq&\frac{1}{p}\| u_k\|^p_{W_\lambda(a)}+\frac{1}{q}\| v_k\|^q_{W_\lambda(b)}-\frac{1}{\alpha}\int_V\left(F_{u_k}(x,u_k,v_k)u_k+F_{v_k}(x,u_k,v_k)v_k+\varrho(|u_k|^p+|v_k|^q)\right)d\mu\nonumber\\
&=&\left(\frac{1}{p}-\frac{1}{\alpha}\right)\| u_k\|^p_{W_\lambda(a)}+\left(\frac{1}{q}-\frac{1}{\alpha}\right)\| v_k\|^q_{W_\lambda(b)}-\frac{\varrho}{\alpha}\| u_k\|^{p}_p-\frac{\varrho}{\alpha}\|v_k\|^{q}_q\nonumber\\
&\geq&\left(\frac{1}{p}-\frac{1}{\alpha}-\frac{\varrho}{\alpha}\right)  \| u_k\|^{p}_{W_\lambda(a)}
+ \left(\frac{1}{q}-\frac{1}{\alpha}-\frac{\varrho}{\alpha}\right)\|v_k\|^{q}_{W_\lambda(b)}.
\end{eqnarray*}
Since $0  <\varrho < \min\{\frac{\alpha-p}{p }, \frac{\alpha-q}{q}\}$, we infer that $\{(u_k,v_k)\}\subseteq \mathcal{N_\lambda}$ is bounded. Therefore, $J_\lambda|_\mathcal{N_\lambda}$ is coercive.
\qed

 \vskip2mm
\noindent
{\bf Lemma 3.7.} {\it If $(F_1)$, $(F_3)$ and $(F_4)$ hold, then there exists some $(u_\lambda,v_\lambda)\in\mathcal{N_\lambda}$ such that $m_\lambda$ can be achieved.}
\vskip2mm
\noindent
{\bf Proof.}\ \ Let $\{(u_k, v_k)\} \subseteq\mathcal{N_\lambda}$ be a minimizing sequence, that is,
$$
\lim\limits_{k\rightarrow +\infty} J_\lambda(u_k,v_k)=m_\lambda.
$$
From Lemma 3.6, we have $\{(u_k,v_k)\}\subseteq \mathcal{N_\lambda}$ is bounded, i.e there exists a positive constant $T$ such that $\|(u_k,v_k)\|_{W_\lambda}\leq T$. Then
\begin{eqnarray}
\label{S16}
T
\geq\| u_k\|_{W_\lambda(a)}
\geq\left(\int_V|u|^p d\mu\right)^{\frac{1}{p}}
\geq\|u_k\|_{\infty,V}.
\end{eqnarray}
Similarly, $\|v_k\|_{\infty,V} \leq T$.
Lemma 2.3 tells us that there exists some $(u_\lambda,v_\lambda) \in W_\lambda$ such that, up to a subsequence,
\begin{eqnarray}
\label{eq3}
 \begin{cases}
  (u_k,v_k)\rightharpoonup (u_\lambda,v_\lambda)\;\;\;\;\;\;\;\;\;\;\;\;\;\;\;\;\;\;\;\;\;\;\mbox{in}\; W_\lambda,\\
   u_k(x)\rightarrow u_\lambda(x),v_k(x)\rightarrow v_\lambda(x)\;\;\;\; \forall x\in V,\\
    u_k \rightarrow  u_\lambda \;\;\;\;\;\;\;\;\;\;\;\;\;\;\;\;\;\;\;\;\;\;\;\;\;\;\;\;\;\;\;\;\;\;\;\;\;\;\mbox{in}\; L^{\theta_1}(V), \;\;\forall \theta_1\in[1,+\infty],\\
    v_k \rightarrow v_\lambda   \;\;\;\;\;\;\;\;\;\;\;\;\;\;\;\;\;\;\;\;\;\;\;\;\;\;\;\;\;\;\;\;\;\;\;\;\;\;\;\mbox{in}\; L^{\theta_2}(V), \;\;\forall \theta_2\in[1,+\infty],
   \end{cases}
\end{eqnarray}
as $k\rightarrow +\infty$.
By $(F_3)$, (\ref{S16}) and Appendix A.3, there exists positive functions $C_3(x)=\max\{1+\frac{1}{p}+\frac{q-1}{q}, 1+\frac{1}{q}+\frac{p-1}{p}\}\cdot C_1(x) , C_4(x)=\max\{1+\frac{1}{r_1}+\frac{r_2-1}{r_2}, 1+\frac{1}{r_2}+\frac{r_1-1}{r_1}\}\cdot C_2(x)$ such that
\begin{eqnarray*}
\label{Z1}
&&|F_{u_k}(x,u_k,v_k)u_k+F_{v_k}(x,u_k,v_k)v_k |\nonumber\\
&\leq & C_3(x)(|u_k|^p+|v_k|^q)+C_4(x)(|u_k|^{r_1}+|v_k|^{r_2}) \nonumber\\
&\leq&   C_3(x)(\|u_k\|_{\infty,V}^p+\|v_k\|_{\infty,V}^q)+C_4(x)(\|u_k\|_{\infty,V}^{r_1}+\|v_k\|_{\infty,V}^{r_2})\nonumber\\
&\leq&   C_3(x)(T^p+T^q)+C_4(x)(T^{r_1}+T^{r_2})\nonumber\\
&\leq&   \max\{T^p+T^q, N^{r_1}+T^{r_2}\}(C_3(x)+C_4(x))\nonumber\\
&:=&g_1(x).
\end{eqnarray*}
Note that $C_1,C_2\in L^1(V)$. Then we can obtain $g_1(x)\in L^1(V)$.
Thus from Lebesgue dominated convergence theorem and (\ref{eq3}), we have
\begin{eqnarray}
\label{Z15}
\lim\limits_{k\rightarrow +\infty} \int_V [F_{u_k}(x,u_k,v_k)u_k+F_{v_k}(x,u_k,v_k)v_k]d\mu= \int_V[F _{u_\lambda}(x,u_\lambda,v_\lambda)u_\lambda+F_{v_\lambda}(x,u_\lambda,v_\lambda)v_\lambda]d\mu.
\end{eqnarray}
Simiarly, by $(F_3)$, (\ref{S16}), (\ref{eq3}), Appendix A.2 and Lebesgue dominated convergence theorem, we also have
\begin{eqnarray}
\label{Z4}
\lim\limits_{k\rightarrow +\infty} \int_V F(x,u_k,v_k)d\mu= \int_V F (x,u_\lambda,v_\lambda)d\mu.
\end{eqnarray}

Since $(u_k, v_k) \in \mathcal{N_\lambda}$, we have
\begin{eqnarray}
\label{Z14}
\| u_k\|^p_{W_\lambda(a)}+\| v_k\|^q_{W_\lambda(b)} =\int_V [F_{u_k}(x,u_k,v_k)u_k+F_{v_k}(x,u_k,v_k)v_k]d\mu.
\end{eqnarray}
By the weak lower continuity of the norm, (\ref{Z15}) and (\ref{Z14}), we have
\begin{eqnarray}
\label{t8}
\| u_\lambda\|^p_{W_\lambda(a)}+ \| v_\lambda\|^q_{W_\lambda(b)} &\leq& \liminf\limits_{k\rightarrow+\infty}\left(\| u_k\|^p_{W_\lambda(a)}+ \| v_k\|^q_{W_\lambda(b)}\right)\nonumber\\
&=&\liminf\limits_{k\rightarrow +\infty}\left(\int_VF_{u_k}(x,u_k,v_k)u_kd\mu+\int_VF_{v_k}(x,u_k,v_k)v_kd\mu\right)\nonumber\\
&=&\int_VF_{u_\lambda}(x,u_\lambda,v_\lambda)u_\lambda d\mu+\int_VF_{v_\lambda}(x,u_\lambda,v_\lambda)v_\lambda d\mu.
\end{eqnarray}
Next, we show that $(u_\lambda, v_\lambda) \in \mathcal{N_\lambda}$ and $(u_\lambda, v_\lambda) \neq (0,0)$.
\par
If $(u_\lambda, v_\lambda) =(0,0)$, then from (\ref{Z15}) and (\ref{Z14}), we have
\begin{eqnarray*}
\liminf\limits_{k\rightarrow +\infty}\left(\| u_k\|^p_{W_\lambda(a)}+\| v_k\|^q_{W_\lambda(b)}\right)&=&\liminf\limits_{k\rightarrow +\infty}\int_V [F_{u_k}(x,u_k,v_k)u_k+F_{v_k}(x,u_k,v_k)v_k]d\mu\\
&= &\int_V[F _{u_\lambda}(x,u_\lambda,v_\lambda)u_\lambda+F_{v_\lambda}(x,u_\lambda,v_\lambda)v_\lambda]d\mu\\
&=&0.
\end{eqnarray*}
Thus, $(u_k,v_k)\rightarrow (0,0)$ in $W_\lambda$ and then by the continuity of $J_\lambda$, (\ref{Z4}) and $(F_1)$ we get $J_\lambda(u_k, v_k)\rightarrow m_\lambda=0.$ It is a contradiction to Lemma 3.5. Therefore, $(u_\lambda, v_\lambda) \neq (0,0)$.
\par
Moreover, from (\ref{t8}), we suppose that
\begin{eqnarray}
\label{b3}
\| u_\lambda\|^p_{W_\lambda(a)}+ \| v_\lambda\|^q_{W_\lambda(b)}<\int_V[F_{u_\lambda}(x,u_\lambda,v_\lambda)u_\lambda+ F_{v_\lambda}(x,u_\lambda,v_\lambda)v_\lambda ]d\mu,
\end{eqnarray}
which implies that $g'(1)<0$. By Lemma 3.3, we know that for $(u_\lambda, v_\lambda)\in W_\lambda \setminus \{(0,0)\}$, there is a unique $t_\lambda$ such that $(t_\lambda u_\lambda, t_\lambda v_\lambda) \in \mathcal{N_\lambda}$ and $g'(t_\lambda)=0$. By Lemma 3.4 and $(F_4)$, we get
\begin{eqnarray*}
&&g'(t_\lambda)
=  t_\lambda^{p-1}\| u_\lambda\|^p_{W_\lambda(a)}+ t_\lambda^{q-1}\|v_\lambda \|^q_{W_\lambda(b)}-\int_V [F_{t_\lambda u_\lambda}(x,t_\lambda u_\lambda,t_\lambda v_\lambda)u_\lambda+F_{t_\lambda v_\lambda}(x,t_\lambda u_\lambda,t_\lambda v_\lambda)v_\lambda]d\mu
=0\\
&\Rightarrow&t_\lambda^{p-1}\| u_\lambda\|^p_{W_\lambda(a)}+ t_\lambda^{q-1}\|v_\lambda \|^q_{W_\lambda(b)}=\int_V [F_{t_\lambda u_\lambda}(x,t_\lambda u_\lambda,t_\lambda v_\lambda)u_\lambda+F_{t_\lambda v_\lambda}(x,t_\lambda u_\lambda,t_\lambda v_\lambda)v_\lambda]d\mu,
\end{eqnarray*}
and it follows from (\ref{b3}), $(F_1)$ and $(F_4)$ that
\begin{eqnarray*}
g''(t_\lambda)
&=  &(p-1)t_\lambda^{p-2}\| u_\lambda\|^p_{W_\lambda(a)}+(q-1)t_\lambda^{q-2}\|v_\lambda \|^q_{W_\lambda(b)}\\
&&- \int_V [F_{t_\lambda u_\lambda,t_\lambda u_\lambda}(x,t_\lambda u_\lambda,t_\lambda v_\lambda)u_\lambda^2+F_{t_\lambda v_\lambda,t_\lambda v_\lambda}(x,t_\lambda u_\lambda,t_\lambda v_\lambda)v_\lambda^2\\
&&-\int_V [F_{t_\lambda u_\lambda,t_\lambda v_\lambda}(x,t_\lambda u_\lambda,t_\lambda v_\lambda)u_\lambda v_\lambda+F_{t_\lambda v_\lambda,t_\lambda u_\lambda}(x,t_\lambda u_\lambda,t_\lambda v_\lambda)v_\lambda u_\lambda]d\mu\\
&  <  & (p-1)t_\lambda^{p-2}\| u_\lambda\|^p_{W_\lambda(a)}+(q-1)t_\lambda^{q-2}\|v_\lambda \|^q_{W_\lambda(b)}\\
&&- (\beta-1)\int_V [\frac{1}{t_\lambda}F_{t_\lambda u_\lambda}(x,t_\lambda u_\lambda,t_\lambda v_\lambda)u_\lambda
+\frac{1}{t_\lambda}F_{t_\lambda v_\lambda}(x,t_\lambda u_\lambda,t_\lambda v_\lambda)v_\lambda]d\mu\\
&&-\int_V [F_{t_\lambda u_\lambda,t_\lambda v_\lambda}(x,t_\lambda u_\lambda,t_\lambda v_\lambda)u_\lambda v_\lambda+F_{t_\lambda v_\lambda,t_\lambda u_\lambda}(x,t_\lambda u_\lambda,t_\lambda v_\lambda)u_\lambda v_\lambda ]d\mu.\\
& <  & 0.
\end{eqnarray*}

Thus we obtain that $t_\lambda$ is the unique maximizer of $g(t)$. Then together with $g'(1)<0$ and $g(0)=0$ we can infer that there exists $t_\lambda \in (0,1)$ such that $(t_\lambda u_\lambda,t_\lambda v_\lambda)\in\mathcal{N_\lambda}$.

\par
When $\beta =p$, by Appendix A.4, (\ref{Z15}), (\ref{Z4}) and the weak lower continuity of the norm, we get
\begin{eqnarray}
\label{t9}
&&m_\lambda\nonumber\\
&\leq& J_\lambda(t_\lambda u_\lambda, t_\lambda v_\lambda)\nonumber\\
&=&\frac{1}{p}\| t_\lambda u_\lambda\|^p_{W_\lambda(a)}+\frac{1}{q}\| t_\lambda v_\lambda\|^q_{W_\lambda(b)} -\int_V F(x,t_\lambda u_\lambda,t_\lambda v_\lambda)d\mu\nonumber\\
&=&\frac{1}{p}\left[\int_V F_{t_\lambda u_\lambda}(x,t_\lambda u_\lambda ,t_\lambda v_\lambda )t_\lambda u_\lambda d\mu+\int_V F_{t_\lambda v_\lambda}(x,t_\lambda u_\lambda ,t_\lambda v_\lambda )t_\lambda v_\lambda d\mu -\| t_\lambda v_\lambda\|^q_{W_\lambda(b)}\right]\nonumber\\
&&+\frac{1}{q}\| t_\lambda v_\lambda\|^q_{W_\lambda(b)}-\int_V F(x,t_\lambda u_\lambda ,t_\lambda v_\lambda )d\mu\nonumber\\
&=&\int_V\left[\frac{1}{p} F_{t_\lambda u_\lambda}(x, t_\lambda u_\lambda ,t_\lambda v_\lambda )t_\lambda u_\lambda +\frac{1}{p} F_{t_\lambda v_\lambda}(x,t_\lambda u_\lambda ,t_\lambda v_\lambda )t_\lambda v_\lambda - F(x,t_\lambda u_\lambda ,t_\lambda v_\lambda )\right]d\mu\nonumber\\
&&+\left(\frac{1}{q}-\frac{1}{p}\right)\| t_\lambda v_\lambda\|^q_{W_\lambda(b)}\nonumber\\
&<&\int_V\left[\frac{1}{p} F_{u_\lambda}(x,u_\lambda ,v_\lambda )u_\lambda +\frac{1}{p} F_{v_\lambda}(x, u_\lambda ,v_\lambda ) v_\lambda - F(x, u_\lambda , v_\lambda )\right]d\mu\\
&&+\left(\frac{1}{q}-\frac{1}{p}\right)\|  v_\lambda\|^q_{W_\lambda(b)}\nonumber\\
&\leq&\liminf\limits_{k\rightarrow +\infty}\left[\int_V \left(\frac{1}{p}F_{u_k}(x, u_k ,v_k ) u_k+\frac{1}{p}F_{v_k}(x, u_k , v_k )v_k -F(x, u_k, v_k)\right)d\mu\right.\nonumber\\
&&\left.+\left(\frac{1}{q}-\frac{1}{p}\right)\|v_k\|^q_{W_\lambda(b)}\right]\nonumber\\
&=&\liminf\limits_{k\rightarrow +\infty}\left[\frac{1}{p}\|u_k\|^p_{W_\lambda(a)}+\frac{1}{p}\|v_k\|^q_{W_\lambda(b)}+\left(\frac{1}{q}-\frac{1}{p}\right)\|v_k\|^q_{W_\lambda(b)} -\int_V F(x,u_k,v_k)d\mu\right]\nonumber\\
&=&m_\lambda\nonumber
\end{eqnarray}
which is a contradiction.
When $\beta =q$, we also have the similar contradiction. Hence (\ref{b3}) does not hold.

Therefore,
$$
\int_VF_{u_\lambda}(x,u_\lambda,v_\lambda)u_\lambda d\mu+\int_VF_{v_\lambda}(x,u_\lambda,v_\lambda)v_\lambda d\mu = \| u_\lambda\|^p_{W_\lambda(a)}+ \| v_\lambda\|^q_{W_\lambda(b)},
$$
and then $( u_\lambda, v_\lambda) \in \mathcal{N}_\lambda$. Thus $m_\lambda$ is achieved by $(u_\lambda,v_\lambda)$.
\qed

 \vskip2mm
 \noindent
{\bf Lemma 3.8.} {\it If $(F_3)$ holds, there exists a constant $\xi>0$ which is independent of $
\lambda$, such that $\|(u,v)\|_{W_\lambda} \geq \xi$ for all $(u,v) \in \mathcal{N}_\lambda$.}

\vskip2mm
\noindent
{\bf Proof.}\ \  Since $(u,v)\in\mathcal{N}_\lambda$, by Lemma 2.2 and Appendix A.3, we have
\begin{eqnarray}
\label{S10}
0&=& \langle J'_\lambda(u ,v ),(u ,v )\rangle\nonumber\\
&=&\| u\|^p_{W_\lambda(a)}+\| v\|^q_{W_\lambda(b)}-\int_VF_u(x,u,v)ud\mu-\int_VF_v(x,u,v)ud\mu\nonumber\\
&\geq&\| u\|^p_{W_\lambda(a)}+\| v\|^q_{W_\lambda(b)}- \int_V[ C_3(x)(|u |^p+|v |^q)+C_4(x)(|u |^{r_1}+|v |^{r_2})]d\mu\nonumber\\
&=&\| u\|^p_{W_\lambda(a)}+\| v\|^q_{W_\lambda(b)}-\|C_3\|_{\infty,V}\|u\|^p_{p}-\|C_3\|_{\infty,V}\|v\|^q_{q}-\|C_4\|_{\infty,V}\|u\|^{r_1}_{r_1}-\|C_4\|_{\infty,V}\|v\|^{r_2}_{r_2}\nonumber\\
&\geq&\| u\|^p_{W_\lambda(a)}+\| v\|^q_{W_\lambda(b)}\nonumber\\
&&-\|C_3\|_{\infty,V}\|u\|^p_{W_\lambda(a)}-\|C_3\|_{\infty,V}\|v\|^q_{W_\lambda(b)}-\|C_4\|_{\infty,V} K_{(p,r_1)}^{r_1}\|u\|^{r_1}_{W_\lambda(a)}-\|C_4\|_{\infty,V} K_{(q,r_2)}^{r_2}\|v\|^{r_2}_{W_\lambda(b)}\nonumber\\
&=&(1-\|C_3\|_{\infty,V})\|u\|^p_{W_\lambda(a)}+(1-\|C_3\|_{\infty,V})\|v\|^q_{W_\lambda(b)}\nonumber\\
&&-\|C_4\|_{\infty,V} K_{(p,r_1)}^{r_1}\|u\|^{r_1}_{W_\lambda(a)}-\|C_4\|_{\infty,V} K_{(q,r_2)}^{r_2}\|v\|^{r_2}_{W_\lambda(b)}.
\end{eqnarray}
\par
We will discuss the above formula in categories.\\
Case $(1):$ assume that $p>q$, $\|u\|_{W_\lambda(a)}\geq 1$ and $\|v\|_{W_\lambda(b)}\geq 1$. Then $\|(u,v)\|_{W_\lambda}\geq 1$ and by (\ref{S10}),
\begin{eqnarray*}
&&0\\
&\geq&(1-\|C_3\|_{\infty,V})\|u\|^p_{W_\lambda(a)}+(1-\|C_3\|_{\infty,V})\|v\|^q_{W_\lambda(b)}-\|C_4\|_{\infty,V} K_{(p,r_1)}^{r_1}\|u\|^{r_1}_{W_\lambda(a)}-\|C_4\|_{\infty,V} K_{(q,r_2)}^{r_2}\|v\|^{r_2}_{W_\lambda(b)}\\
&\geq&(1-\|C_3\|_{\infty,V})(\|u\|^q_{W_\lambda(a)}+\|v\|^q_{W_\lambda(b)})-\max\{\|C_4\|_{\infty,V} K_{(p,r_1)}^{r_1}, \|C_4\|_{\infty,V} K_{(q,r_2)}^{r_2}\}(\|u\|^{r_1}_{W_\lambda(a)}+\|v\|^{r_2}_{W_\lambda(b)})\\
&\geq&\frac{ 1-\|C_3\|_{\infty,V} }{2^{q-1}}(\|u\|_{W_\lambda(a)}+\|v\|_{W_\lambda(b)})^q-\max\{\|C_4\|_{\infty,V} K_{(p,r_1)}^{r_1}, \|C_4\|_{\infty,V} K_{(q,r_2)}^{r_2}\}(\|(u,v)\|^{r_1}_{W_\lambda}+\|(u,v)\|^{r_2}_{W_\lambda})\\
&\geq&\frac{ 1-\|C_3\|_{\infty,V}}{2^{q-1}}\|(u,v)\|_{W_\lambda}^q-2\max\{\|C_4\|_{\infty,V} K_{(p,r_1)}^{r_1}, \|C_4\|_{\infty,V} K_{(q,r_2)}^{r_2}\}\|(u,v)\|_{W_\lambda}^{\max\{r_1,r_2\}}.
\end{eqnarray*}
Thus we have
\begin{eqnarray*}
\label{S2}
\|(u,v)\|_{W_\lambda} \geq \left(\frac{ 1-\|C_3\|_{\infty,V} }{2^{q}\max\{\|C_4\|_{\infty,V} K_{(p,r_1)}^{r_1}, \|C_4\|_{\infty,V} K_{(q,r_2)}^{r_2}\}}\right)^\frac{1}{\max\{r_1,r_2\}-q}.
\end{eqnarray*}
Case $(2):$ assume that $p>q$, $\|u\|_{W_\lambda(a)}\geq 1$ and $\|v\|_{W_\lambda(b)}< 1$. Then $\|(u,v)\|_{W_\lambda}\geq 1$ and by (\ref{S10}),
\begin{eqnarray*}
&&0\\
&\geq&(1-\|C_3\|_{\infty,V})\|u\|^p_{W_\lambda(a)}+(1-\|C_3\|_{\infty,V})\|v\|^q_{W_\lambda(b)}-\|C_4\|_{\infty,V} K_{(p,r_1)}^{r_1}\|u\|^{r_1}_{W_\lambda(a)}-\|C_4\|_{\infty,V} K_{(q,r_2)}^{r_2}\|v\|^{r_2}_{W_\lambda(b)}\nonumber\\
&\geq&(1-\|C_3\|_{\infty,V})\|u\|^p_{W_\lambda(a)}-\|C_4\|_{\infty,V}  K_{(p,r_1)}^{r_1}\|u\|^{r_1}_{W_\lambda(a)}-\|C_4\|_{\infty,V} K_{(q,r_2)}^{r_2}\|v\|^{r_2}_{W_\lambda(b)}.
\end{eqnarray*}
Thus we have
\begin{eqnarray*}
\label{S3}
1-\|C_3\|_{\infty,V}  &\leq& \|C_4\|_{\infty,V} K_{(p,r_1)}^{r_1}\|u\|^{r_1-p}_{W_\lambda(a)}+\|C_4\|_{\infty,V} K_{(q,r_2)}^{r_2}\frac{\|v\|^{r_2}_{W_\lambda(b)}}{\|u\|^p_{W_\lambda(a)}}\nonumber\\
  & \leq & \|C_4\|_{\infty,V} K_{(p,r_1)}^{r_1}\|u\|^{r_1-p}_{W_\lambda(a)}+\|C_4\|_{\infty,V} K_{(q,r_2)}^{r_2}\|v\|^{r_2}_{W_\lambda(b)}\nonumber\\
  & \leq &\|C_4\|_{\infty,V} K_{(p,r_1)}^{r_1}\|(u,v)\|_{W_\lambda}^{r_1-p}+\|C_4\|_{\infty,V} K_{(q,r_2)}^{r_2}\|(u,v)\|_{W_\lambda}^{r_2}\nonumber\\
  & \leq &2\max\{\|C_4\|_{\infty,V} K_{(p,r_1)}^{r_1}, \|C_4\|_{\infty,V} K_{(q,r_2)}^{r_2}\} \|(u,v)\|_{W_\lambda}^{\max\{r_1-p,r_2\}}\nonumber\\
\Rightarrow  \|(u,v)\|_{W_\lambda} &\geq& \left(\frac{ 1-\|C_3\|_{\infty,V}  }{2\max\{\|C_4\|_{\infty,V} K_{(p,r_1)}^{r_1}, \|C_4\|_{\infty,V} K_{(q,r_2)}^{r_2}\}}\right)^\frac{1}{\max\{r_1-p,r_2\}}.
\end{eqnarray*}
Case $(3):$ assume that $p>q$, $\|u\|_{W_\lambda(a)}< 1$ and $\|v\|_{W_\lambda(b)}\geq1$. Then $\|(u,v)\|_{W_\lambda}\geq 1$ and similar to case $(2)$, we have
\begin{eqnarray*}
\label{S5}
\|(u,v)\|_{W_\lambda} \geq \left(\frac{ 1-\|C_3\|_{\infty,V} }{2\max\{\|C_4\|_{\infty,V} K_{(p,r_1)}^{r_1}, \|C_4\|_{\infty,V} K_{(q,r_2)}^{r_2}\}}\right)^\frac{1}{\max\{r_1,r_2-q\}}.
\end{eqnarray*}
Case $(4):$ assume that $p>q$, $\|u\|_{W_\lambda(a)}< 1$ and $\|v\|_{W_\lambda(b)}< 1$, by (\ref{S10}),
\begin{eqnarray*}
&&0\\
&\geq&(1-\|C_3\|_{\infty,V} )(\|u\|^p_{W_\lambda(a)}+\|v\|^p_{W_\lambda(b)})-\max\{\|C_4\|_{\infty,V} K_{(p,r_1)}^{r_1}, \|C_4\|_{\infty,V} K_{(q,r_2)}^{r_2}\}(\|u\|^{r_1}_{W_\lambda(a)}+\|v\|^{r_2}_{W_\lambda(b)})\\
&\geq&\frac{ 1-\|C_3\|_{\infty,V}  }{2^{p-1}}\|(u,v)\|_{W_\lambda}^p-\max\{\|C_4\|_{\infty,V} K_{(p,r_1)}^{r_1}, \|C_4\|_{\infty,V} K_{(q,r_2)}^{r_2}\}(\|(u,v)\|^{r_1}_{W_\lambda}+\|(u,v)\|^{r_2}_{W_\lambda}).
\end{eqnarray*}
Thus we have
\begin{eqnarray}
\label{S6}
\max\{\|C_4\|_{\infty,V} K_{(p,r_1)}^{r_1}, \|C_4\|_{\infty,V} K_{(q,r_2)}^{r_2}\}(\|(u,v)\|^{r_1-p}_{W_\lambda}+\|(u,v)\|^{r_2-p}_{W_\lambda}) \geq \frac{ 1-\|C_3\|_{\infty,V}  }{2^{p-1}}.
\end{eqnarray}
When $\|(u,v)\|_{W_\lambda}\geq 1$, from (\ref{S6}), we have
\begin{eqnarray*}
\label{S7}
&&2\max\{\|C_4\|_{\infty,V} K_{(p,r_1)}^{r_1}, \|C_4\|_{\infty,V} K_{(q,r_2)}^{r_2}\}\|(u,v)\|^{\max\{r_1-p,r_2-p\}}\geq \frac{ 1-\|C_3\|_{\infty,V}  }{2^{p-1}}\nonumber\\
&\Rightarrow&\|(u,v)\|_{W_\lambda} \geq \left(\frac{ 1-\|C_3\|_{\infty,V}  }{2^{p}\max\{\|C_4\|_{\infty,V} K_{(p,r_1)}^{r_1}, \|C_4\|_{\infty,V} K_{(q,r_2)}^{r_2}\}}\right)^\frac{1}{\max\{r_1,r_2\}-p}.
\end{eqnarray*}
When $\|(u,v)\|_{W_\lambda}< 1$, from (\ref{S6}), we have
\begin{eqnarray*}
\label{S8}
&&2\max\{\|C_4\|_{\infty,V} K_{(p,r_1)}^{r_1}, \|C_4\|_{\infty,V} K_{(q,r_2)}^{r_2}\}\|(u,v)\|^{\min\{r_1-p,r_2-p\}}\geq \frac{ 1-\|C_3\|_{\infty,V}  }{2^{p-1}}\nonumber\\
&\Rightarrow&\|(u,v)\|_{W_\lambda} \geq \left(\frac{1-\|C_3\|_{\infty,V} }{2^{p}\max\{\|C_4\|_{\infty,V} K_{(p,r_1)}^{r_1}, \|C_4\|_{\infty,V} K_{(q,r_2)}^{r_2}\}}\right)^\frac{1}{\min\{r_1,r_2\}-p}.
\end{eqnarray*}
For the case $p<q$, using the same method, we can obtain that\\
Case $(5):$ assume that $\|u\|_{W_\lambda(a)}\geq 1$ and $\|v\|_{W_\lambda(b)}\geq1$, $$\|(u,v)\|_{W_\lambda}\geq\left(\frac{ 1-\|C_3\|_{\infty,V}  }{2^{p}\max\{\|C_4\|_{\infty,V} K_{(p,r_1)}^{r_1}, \|C_4\|_{\infty,V} K_{(q,r_2)}^{r_2}\}}\right)^\frac{1}{\max\{r_1,r_2\}-p};$$
Case $(6):$ assume that $\|u\|_{W_\lambda(a)}\geq 1$ and $\|v\|_{W_\lambda(b)}<1$,
$$\|(u,v)\|_{W_\lambda}\geq\left(\frac{ 1-\|C_3\|_{\infty,V}  }{2\max\{\|C_4\|_{\infty,V} K_{(p,r_1)}^{r_1}, \|C_4\|_{\infty,V} K_{(q,r_2)}^{r_2}\}}\right)^\frac{1}{\max\{r_1-p,r_2\}};$$
Case $(7):$ assume that $\|u\|_{W_\lambda(a)}< 1$ and $\|v\|_{W_\lambda(b)}\geq1$,
$$\|(u,v)\|_{W_\lambda} \geq \left(\frac{ 1-\|C_3\|_{\infty,V}  }{2\max\{\|C_4\|_{\infty,V} K_{(p,r_1)}^{r_1}, \|C_4\|_{\infty,V} K_{(q,r_2)}^{r_2}\}}\right)^\frac{1}{\max\{r_1,r_2-q\}};$$
Case $(8):$ assume that $\|u\|_{W_\lambda(a)}< 1$ and $\|v\|_{W_\lambda(b)}<1$,
\begin{eqnarray*}
\|(u,v)\|_{W_\lambda} \geq \left(\frac{ 1-\|C_3\|_{\infty,V}  }{2^{q}\max\{\|C_4\|_{\infty,V} K_{(p,r_1)}^{r_1}, \|C_4\|_{\infty,V} K_{(q,r_2)}^{r_2}\}}\right)^\frac{1}{\max\{r_1,r_2\}-q},\;\;\;\mbox{for}\;\|(u,v)\|_{W_\lambda}\geq 1.
\end{eqnarray*}
\begin{eqnarray*}
\|(u,v)\|_{W_\lambda} \geq \left(\frac{ 1-\|C_3\|_{\infty,V}  }{2^{q}\max\{\|C_4\|_{\infty,V} K_{(p,r_1)}^{r_1}, \|C_4\|_{\infty,V} K_{(q,r_2)}^{r_2}\}}\right)^\frac{1}{\min\{r_1,r_2\}-q},\;\;\;\mbox{for}\;\|(u,v)\|_{W_\lambda}< 1.
\end{eqnarray*}
Case $(9):$ assume that $p=q$, from (\ref{S10}) we have
\begin{eqnarray*}
&&0\\
&\geq&(1-\|C_3\|_{\infty,V} )\|u\|^p_{W_\lambda(a)}+(1-\|C_3\|_{\infty,V} )\|v\|^q_{W_\lambda(b)}-\|C_4\|_{\infty,V} K_{(p,r_1)}^{r_1}\|u\|^{r_1}_{W_\lambda(a)}-\|C_4\|_{\infty,V} K_{(q,r_2)}^{r_2}\|v\|^{r_2}_{W_\lambda(b)}\\
&\geq&\frac{1}{2^{p-1}}(1-\|C_3\|_{\infty,V} )\|(u,v)\|^p-\max\{\|C_4\|_{\infty,V} K_{(p,r_1)}^{r_1}, \|C_4\|_{\infty,V} K_{(q,r_2)}^{r_2}\}(\|(u,v)\|^{r_1}_{W_\lambda}+\|(u,v)\|^{r_2}_{W_\lambda}).
\end{eqnarray*}
Thus we get
\begin{eqnarray}
\label{S9}
\max\{\|C_4\|_{\infty,V} K_{(p,r_1)}^{r_1}, \|C_4\|_{\infty,V} K_{(q,r_2)}^{r_2}\}(\|(u,v)\|^{r_1-p}_{W_\lambda}+\|(u,v)\|^{r_2-p}_{W_\lambda})\geq  \frac{1}{2^{p-1}} (1-\|C_3\|_{\infty,V})  .
\end{eqnarray}
When $\|(u,v)\|_{W_\lambda}\geq 1$, from (\ref{S9}), there is
\begin{eqnarray*}
\label{S11}
&&2\max\{\|C_4\|_{\infty,V} K_{(p,r_1)}^{r_1}, \|C_4\|_{\infty,V} K_{(q,r_2)}^{r_2}\}\|(u,v)\|^{\max\{r_1-p,r_2-p\}}\geq \frac{1}{2^{p-1}} (1-\|C_3\|_{\infty,V} ) \nonumber\\
&\Rightarrow&\|(u,v)\|_{W_\lambda} \geq \left(\frac{1-\|C_3\|_{\infty,V} }{2^p\max\{\|C_4\|_{\infty,V} K_{(p,r_1)}^{r_1}, \|C_4\|_{\infty,V} K_{(q,r_2)}^{r_2}\}}\right)^\frac{1}{\max\{r_1,r_2\}-p}.
\end{eqnarray*}
When $\|(u,v)\|_{W_\lambda}< 1$, from (\ref{S9}), there is
\begin{eqnarray*}
\label{S12}
&&2\max\{\|C_4\|_{\infty,V} K_{(p,r_1)}^{r_1}, \|C_4\|_{\infty,V} K_{(q,r_2)}^{r_2}\}\|(u,v)\|^{\min\{r_1-p,r_2-p\}}\geq \frac{1}{2^{p-1}} (1-\|C_3\|_{\infty,V})  \nonumber\\
&\Rightarrow&\|(u,v)\|_{W_\lambda} \geq \left(\frac{1-\|C_3\|_{\infty,V} }{2^p\max\{\|C_4\|_{\infty,V} K_{(p,r_1)}^{r_1}, \|C_4\|_{\infty,V} K_{(q,r_2)}^{r_2}\}}\right)^\frac{1}{\min\{r_1,r_2\}-p}.
\end{eqnarray*}
Let
\begin{eqnarray*}
&&\xi_1= \left(\frac{ 1-\|C_3\|_{\infty,V}  }{2^{\min\{p,q\}}\max\{\|C_4\|_{\infty,V} K_{(p,r_1)}^{r_1}, \|C_4\|_{\infty,V} K_{(q,r_2)}^{r_2}\}}\right)^\frac{1}{\max\{r_1,r_2\}-\min\{p,q\}},\\
&&\xi_2=\left(\frac{ 1-\|C_3\|_{\infty,V}  }{2\max\{\|C_4\|_{\infty,V} K_{(p,r_1)}^{r_1}, \|C_4\|_{\infty,V} K_{(q,r_2)}^{r_2}\}}\right)^\frac{1}{\max\{r_1-p,r_2\}},\\
&&\xi_3=\left(\frac{ 1-\|C_3\|_{\infty,V}  }{2\max\{\|C_4\|_{\infty,V} K_{(p,r_1)}^{r_1}, \|C_4\|_{\infty,V} K_{(q,r_2)}^{r_2}\}}\right)^\frac{1}{\max\{r_1,r_2-q\}},\\
&&\xi_4= \left(\frac{ 1-\|C_3\|_{\infty,V}  }{2^{\max\{p,q\}}\max\{\|C_4\|_{\infty,V} K_{(p,r_1)}^{r_1}, \|C_4\|_{\infty,V} K_{(q,r_2)}^{r_2}\}}\right)^\frac{1}{\max\{r_1,r_2\}-\max\{p,q\}}.
\end{eqnarray*}
All in all, we can obtain that $\|(u,v)\|_{W_\lambda} \geq \xi$, where $\xi=\min\{\xi_1,\xi_2,\xi_3,\xi_4\}$. By Appendix A.1 and noting that $C_1(x) \leq \frac{1}{1+\beta}$, we obtain that $\|C_3\|_{\infty,V}<1$ and then $\xi>0$. Thus the lemma is proved.
\qed

\vskip2mm
 \noindent
{\bf Lemma 3.9.} {\it If $(F_2)$ holds, then for a minimizing sequence $\{(u_k,v_k)\}\subseteq\mathcal{N}_{\lambda}$, there holds
\begin{eqnarray*}
\label{f2}
\lim \limits_{k \rightarrow +\infty}\| (u_k, v_k)\|_{W_\lambda}\leq \left(\max\{2^{p-1},2^{q-1}\}\left(\frac{m_\lambda}{ \frac{1}{\beta}-\frac{1+\varrho}{\alpha}}+1\right)\right)^{\frac{1}{\min\{p,q\}}}.
\end{eqnarray*}
}

\vskip2mm
\noindent
{\bf Proof.}\ \
 Since $\{(u_k,v_k)\}\subseteq\mathcal{N}_{\lambda}$ is a minimizing sequence, we have
\begin{eqnarray}
\label{Z9}
0=\| u_k\|^p_{W_\lambda(a)}+\| v_k\|^q_{W_\lambda(b)}-\int_VF_{u_k}(x,u_k,v_k)u_k d\mu-\int_VF_{v_k}(x,u_k,v_k)u_k d\mu.
\end{eqnarray}
Then by $(F_2)$ and (\ref{Z9}),
\begin{eqnarray*}
\label{S13}
m_\lambda&   =& \lim \limits_{k\rightarrow +\infty} \left(\frac{1}{p}\| u_k\|^p_{W_\lambda(a)}+\frac{1}{q}\| v_k\|^q_{W_\lambda(b)}-\int_V F(x,u_k,v_k)d\mu\right)\nonumber\\
& \geq&   \lim \limits_{k \rightarrow +\infty}\left[  \frac{1}{p}\| u_k\|^p_{W_\lambda(a)}+\frac{1}{q}\| v_k\|^q_{W_\lambda(b)}\nonumber\right.\\
&&\left.  -\frac{1}{\alpha}\int_V(F_{u_k}(x,u_k,v_k)u_k +F_{v_k}(x,u_k,v_k)v_k+\varrho(|u_k|^p+|v_k|^q)) d\mu\right]\nonumber\\
&       =       & \lim \limits_{k \rightarrow +\infty} \left( \frac{1}{p}\| u_k\|^p_{W_\lambda(a)}+ \frac{1}{q} \| v_k\|^q_{W_\lambda(b)}\right.\nonumber\\
& &  \left.  -\frac{1}{\alpha}\| u_k\|^p_{W_\lambda(a)}-\frac{1}{\alpha}\| v_k\|^q_{W_\lambda(b)}-\frac{\varrho}{\alpha}\|u_k\|^p_p-\frac{\varrho}{\alpha}\|v_k\|^q_q\right)\nonumber\\
&   \geq    & \lim \limits_{k \rightarrow +\infty}\left[  \left(\frac{1}{p}-\frac{1}{\alpha}-\frac{\varrho}{\alpha}\right)\| u_k\|^p_{W_\lambda(a)}+\left(\frac{1}{q}-\frac{1}{\alpha}-\frac{\varrho}{\alpha}\right)\| v_k\|^q_{W_\lambda(b)}\right]
\nonumber\\
&   \geq    &   \lim \limits_{k \rightarrow +\infty}\left[ \min\left\{ \frac{1}{p}-\frac{1+\varrho}{\alpha},\frac{1}{q}-\frac{1+\varrho}{\alpha}\right\}\left(\| u_k\|^p_{W_\lambda(a)}+\| v_k\|^q_{W_\lambda(b)}\right)\right]\nonumber\\
&   \geq    &  \lim \limits_{k \rightarrow +\infty}\left[   \min\left\{ \frac{1}{p}-\frac{1+\varrho}{\alpha},\frac{1}{q}-\frac{1+\varrho}{\alpha}\right\}
\left(\frac{\|(u_k,v_k)\|_{W_\lambda}^{\min\{p,q\}}}{\max\{2^{p-1},2^{q-1}\}}-1\right)\right]\;\;(\mbox{see}\; (66)\; \mbox{in}\; \cite{Xie J 2018}).
\end{eqnarray*}
Then we can obtain that
\begin{eqnarray*}
\label{S14}
\lim \limits_{k \rightarrow +\infty}   \|(u_k,v_k)\|_{W_\lambda} \leq \left(\max\{2^{p-1},2^{q-1}\}\left(\frac{m_\lambda}{ \frac{1}{\beta}-\frac{1+\varrho}{\alpha}}+1\right)\right)^{\frac{1}{\min\{p,q\}}}.
\end{eqnarray*}
The proof is complete.
\qed

\vskip2mm
\noindent
{\bf Proof of Theorem 1.1.}\ \ We denote the critical set
$$
K_{J_\lambda}=\{(u,v)\in W_\lambda:J'_\lambda(u,v)=0\}.
$$
Let $k(u,v):=\langle J'_\lambda(u, v),(u, v) \rangle$, that is,
\begin{eqnarray*}
\label{Z16}
k(u,v)= \| u\|^p_{W_\lambda(a)}+ \| v\|^q_{W_\lambda(b)}-\int_VF_{u }(x,u ,v )u d\mu-\int_VF_{v}(x,u,v)v d\mu .
\end{eqnarray*}
Then by $(F_1)$, $k(u,v)\in C^1(W_\lambda,\R)$ and
\begin{eqnarray*}
\label{Z17}
&&\langle k'(u,v),(h_1, h_2)\rangle\nonumber\\
&=& p\int_V(|\nabla u|^{p-2}\Gamma(u,h_1)+(\lambda a+1)|u|^{p-2}uh_1)d\mu
+q\int_V(|\nabla v|^{q-2}\Gamma(v, h_2)+(\lambda b+1)|v|^{q-2}vh_2)d\mu\nonumber\\
&&-\int_V [F_{uu}(x,u,v)uh_1+F_{uv}(x,u,v)uh_2+F_u(x,u,v)h_1]d\mu\nonumber\\
&&-\int_V [F_{vv}(x,u,v)vh_2+F_{v u}(x,u,v)vh_1+F_v(x,u,v)h_2]d\mu,
\end{eqnarray*}
for all $(h_1,h_2)\in W_\lambda$. From Lemma 3.7, we know that
$$
J_\lambda(u_\lambda, v_\lambda)=m_\lambda=\mbox{inf}\{J_\lambda(u, v):k(u, v)=0, (u, v)\in W_\lambda\setminus{(0,0)}\}.
$$
Then by the Lagrange multiplier rule, we can find $\vartheta\in \R$ such that for $ ( u_\lambda, v_\lambda) \in \mathcal{N}_\lambda $, we have
\begin{eqnarray}
\label{Z18}
&&J'_\lambda(u_\lambda, v_\lambda)+\vartheta k'(u_\lambda,v_\lambda) =0\nonumber\\
&\Rightarrow&\langle J'_\lambda(u_\lambda, v_\lambda), (u_\lambda, v_\lambda) \rangle+\vartheta \langle k'(u_\lambda,v_\lambda),(u_\lambda, v_\lambda)\rangle=0\nonumber\\
&\Rightarrow& \vartheta \langle k'(u_\lambda,v_\lambda),(u_\lambda, v_\lambda)\rangle=0.
\end{eqnarray}
Suppose $\vartheta \neq 0$. Then $ \langle k'(u_\lambda,v_\lambda),(u_\lambda, v_\lambda)\rangle=0$. When $\beta=p$, by $(F_4)$, we get
\begin{eqnarray*}
\label{Z19}
&&\langle k'(u_\lambda,v_\lambda),(u_\lambda, v_\lambda)\rangle\nonumber\\
&= &p\|u_\lambda\|^p_{W_\lambda(a)}+q\|v_\lambda\|^q_{W_\lambda(b)}-\int_V [F_{u_\lambda}(x,u_\lambda,v_\lambda)u_\lambda+F_{v_\lambda}(x,u_\lambda,v_\lambda)v_\lambda]d\mu\nonumber\\
&&-\int_V \left[F_{u_\lambda u_\lambda}(x,u_\lambda,v_\lambda)u^2_\lambda+F_{v_\lambda v_\lambda}(x,u_\lambda,v_\lambda)v^2_\lambda+2F_{u_\lambda v_\lambda}(x,u_\lambda,v_\lambda)u_\lambda v_\lambda\right]d\mu\nonumber\\
&=&p\left[\|u_\lambda\|^p_{W_\lambda(a)}+\|v_\lambda\|^q_{W_\lambda(b)}-\int_V \left[F_{u_\lambda}(x,u_\lambda,v_\lambda)u_\lambda+F_{v_\lambda}(x,u_\lambda,v_\lambda)v_\lambda\right]d\mu\right]\nonumber\\
&&-\int_V \left[F_{u_\lambda u_\lambda}(x,u_\lambda,v_\lambda)u^2_\lambda+F_{v_\lambda v_\lambda}(x,u_\lambda,v_\lambda)v^2_\lambda+2F_{u_\lambda v_\lambda}(x,u_\lambda,v_\lambda)u_\lambda v_\lambda\right]d\mu\nonumber\\
&&+(q-p)\|v_\lambda\|^q_{W_\lambda(b)}+(p-1)\int_V [F_{u_\lambda}(x,u_\lambda,v_\lambda)u_\lambda+F_{v_\lambda}(x,u_\lambda,v_\lambda)v_\lambda]d\mu\nonumber\\
&=& (q-p)\|v_\lambda\|^q_{W_\lambda(b)}-\int_V \left[F_{u_\lambda u_\lambda}(x,u_\lambda,v_\lambda)u^2_\lambda+F_{v_\lambda v_\lambda}(x,u_\lambda,v_\lambda)v^2_\lambda\right.\nonumber\\
&&\left.+2F_{u_\lambda v_\lambda}(x,u_\lambda,v_\lambda)u_\lambda v_\lambda\right]d\mu+(p-1)\int_V \left[F_{u_\lambda}(x,u_\lambda,v_\lambda)u_\lambda+F_{v_\lambda}(x,u_\lambda,v_\lambda)v_\lambda\right]d\mu \nonumber\\
&<&0.
\end{eqnarray*}
It is a contradiction. Similarly, when $\beta =q$, we can also have the same contradiction. So $\vartheta = 0$. From (\ref{Z18}), we get
$J'_\lambda(u_\lambda, v_\lambda)=0$ and so $(u_\lambda, v_\lambda)$ is a nontrivial solution of (\ref{eq1}).
Moreover, by Lemma 3.8, the weak lower semi-continuity of the norm $\|\cdot\|_{W_\lambda}$ and Lemma 3.9, it is easy to obtain that
\begin{eqnarray}
\label{**}
\xi \leq \|(u_\lambda,v_\lambda)\|_{W_\lambda} \leq\left(\max\{2^{p-1},2^{q-1}\}\left(\frac{m_\lambda}{ \frac{1}{\beta}-\frac{1+\varrho}{\alpha}}+1\right)\right)^{\frac{1}{\min\{p,q\}}}.
\end{eqnarray}
Thus the proof is finished.\qed

\vskip2mm
{\section{Convergence of the ground state solutions family}}
  \setcounter{equation}{0}
  \par
In this section, we prove that the ground state solutions family $\{(u_\lambda,v_\lambda)\}$ of (\ref{eq1}) converge to a ground state solution of (\ref{e5}) as $\lambda \rightarrow +\infty$, which imply Theorem 1.2.

 \vskip2mm
 \noindent
{\bf Lemma 4.1.}  {\it If $(F_1)$ and $(F_2)$ hold, then $m_\lambda\rightarrow m_\Omega$ as $\lambda \rightarrow\infty$, where $m_\Omega$ is defined by (\ref{e8}).}

 \vskip0mm
 \noindent
{\bf Proof.} Since $\mathcal{N}_\Omega  \subset  \mathcal{N}_\lambda $(see Appendix A.5), we obviously have that $m_\lambda \leq m_\Omega$ for any $\lambda > 0$. Take a sequence $\lambda_k \rightarrow \infty$ as $k\rightarrow\infty$ such that
\begin{eqnarray}
\label{Z10}
\lim\limits_{k \rightarrow \infty}m_{\lambda_k}=M\leq m_\Omega <+\infty,
\end{eqnarray}
where $m_{\lambda_k}$ is the ground state of (\ref{eq1}) with $\lambda=\lambda_k$ and $(u_{\lambda_k},v_{\lambda_k}) \in \mathcal{N}_{\lambda_k}$ is the corresponding ground state solution.
Similar to the proof of Lemma 3.3, it can be  obtain that $\mathcal{N}_\Omega $ is non-empty, that is, there exists a $(u_*,v_*)\in N_\Omega $. Then by the definition of $J_\Omega$, we obtain that $ m_\Omega=\inf_{(u,v)\in \mathcal{N}_\Omega}J_\Omega(u,v) \le J_\Omega(u_*,v_*) <+\infty$. Lemma 3.5 tells us that $m_\Omega \geq M\ge\eta> 0$.
\par
From (\ref{**}) and (\ref{Z10}), we get
\begin{eqnarray}
\label{*}\|(u_{\lambda_k},v_{\lambda_k})\|_{W_{\lambda_k}} \leq\left(\max\{2^{p-1},2^{q-1}\}\left(\frac{m_{\Omega}}{ \frac{1}{\beta}-\frac{1+\varrho}{\alpha}}+1\right)\right)^{\frac{1}{\min\{p,q\}}}:=L.
\end{eqnarray}
Besides by the definitions of $W$ and $W_\lambda$, it is easy to see that $\{(u_{\lambda_k},v_{\lambda_k})\}\subset W$ and so $\{(u_{\lambda_k},v_{\lambda_k})\} $ is uniformly bounded in $W$. So there exists $(u_0,v_0)\subset W$ such that
\begin{eqnarray}
\label{t5}
 \begin{cases}
  (u_{\lambda_k},v_{\lambda_k})\rightharpoonup (u_0,v_0),\;\;\;\;\;\;\;\;\;\;\;\;\;\;\;\;\;\;\;\;\;\;\;\;\;\;\;\;\;\;\;\;\;\;\;\; \mbox{in}\; W,\\
   u_{\lambda_k}(x)\rightarrow u_0(x),v_{\lambda_k}(x)\rightarrow v_0(x)\;\;\;\;\;\;\;\;\;\;\;\;\;\;\;\;\; \forall x\in V.
   \end{cases}
\end{eqnarray}
Here we have used the proof of Lemma 2.3 with replacing $W_\lambda(a), W_\lambda(b)$ and $W_\lambda$ with $W^{1,p}(V), W^{1,q}(V)$ and $W$, respectively, which is easily verified.
\par
We claim that $u_0|_{\Omega^c_a}\equiv0$ and $v_0|_{\Omega^c_b}\equiv0$. Otherwise, without loss of generality, we assume that there exists a vertex $x_0 \in \Omega^c_a$ such that $u_0(x_0)\neq 0$. Since $(u_{\lambda_k},v_{\lambda_k}) \in \mathcal{N}_{\lambda_k}$, we have
\begin{eqnarray*}
&&J_{\lambda_k}(u_{\lambda_k}, v_{\lambda_k})\\
&   =&  J_{\lambda_k}(u_{\lambda_k}, v_{\lambda_k})-\frac{1}{\alpha}\langle J'_{\lambda_k}(u_{\lambda_k},v_{\lambda_k}),(u_{\lambda_k},v_{\lambda_k})\rangle\\
&       =       & \left(\frac{1}{p}-\frac{1}{\alpha}\right)\| u_{\lambda_k}\|^p_{W_\lambda(a)}+\left(\frac{1}{q}-\frac{1}{\alpha}\right)\| v_{\lambda_k}\|^q_{W_\lambda(b)} \\
& &   +\frac{1}{\alpha}\int_VF_{u_{\lambda_k}}(x,u_{\lambda_k},v_{\lambda_k})u_{\lambda_k}d\mu+\frac{1}{\alpha}\int_VF_{v_{\lambda_k}}(x,u_{\lambda_k},v_{\lambda_k})v_{\lambda_k}d\mu-\int_V F(x,u_{\lambda_k},v_{\lambda_k})d\mu \\
&   \geq    &  \left(\frac{1}{p}-\frac{1}{\alpha}\right)\| u_{\lambda_k}\|^p_{W_\lambda(a)}+\left(\frac{1}{q}-\frac{1}{\alpha}\right)\| v_{\lambda_k}\|^q_{W_\lambda(b)}
-\frac{\varrho}{\alpha}\int_V(|u_{\lambda_k}|^{p}+|v_{\lambda_k}|^{q})d\mu \\
&   \geq   &   \left(\frac{1}{\beta}-\frac{1}{\alpha}-\frac{\varrho}{\alpha}\right) (\| u_{\lambda_k}\|^p_{W_\lambda(a)}+\| v_{\lambda_k}\|^q_{W_\lambda(b)})\\
&     \geq   &  \left(\frac{1}{\beta}-\frac{1}{\alpha}-\frac{\varrho}{\alpha}\right) \int_V(|\nabla u_{\lambda_k}|^p+(\lambda_k a(x)+1)|u_{\lambda_k}|^p)d\mu\\
&   \geq    &  \left(\frac{1}{\beta}-\frac{1}{\alpha}-\frac{\varrho}{\alpha}\right) \lambda_k\int_Va(x)|u_{\lambda_k}|^pd\mu\\
&   \geq    &\left(\frac{1}{\beta}-\frac{1}{\alpha}-\frac{\varrho}{\alpha}\right)  \lambda_k a(x_0)|u_{\lambda_k}(x_0)|^p\mu(x_0) .
\end{eqnarray*}
Since $x_0 \notin \Omega_a, a(x_0) >0$, $\mu(x_0) \geq \mu_{\min}>0$ and $u_{\lambda_k}(x_0) \rightarrow u_0(x_0)\neq 0 $, then by $\lambda_k \rightarrow +\infty$, we get that $\lim\limits_{k \rightarrow \infty}J_{\lambda_k}(u_{\lambda_k}, v_{\lambda_k})= +\infty$, which is a contraction to the fact that $m_{\lambda_k}\leq m_\Omega$. So $u_0|_{\Omega^c_a}\equiv0$. Similarly, we can also obtain that $v_0|_{\Omega^c_b}\equiv0$.
\par
Similar to (\ref{Z4}), by (\ref{t5}), Appendix A.2, Appendix A.3 and Lebesgue dominated convergence theorem, for all $t>0$, we have
\begin{eqnarray}
\label{t7}
\lim\limits_{k \rightarrow \infty}\int_VF(x,tu_{\lambda_k},tv_{\lambda_k})d\mu=\int_VF(x,tu_0,tv_0)d\mu
\end{eqnarray}
\begin{eqnarray}
\label{a3}
\lim\limits_{k \rightarrow \infty}\int_V(F_{tu_{\lambda_k}}(x,tu_{\lambda_k},tv_{\lambda_k})tu_{\lambda_k}+F_{tv_{\lambda_k}}(x,tu_{\lambda_k},tv_{\lambda_k})tv_{\lambda_k})d\mu=\int_V \left(F_{tu_{0}}(x,tu_{0},tv_0)tu_{0}+F_{tv_0}(x,tu_0,tv_0)tv_0\right)d\mu.
\end{eqnarray}
\par
Since $(u_0, v_0)\in W$ and $u_0|_{\Omega^c_a}\equiv0$ and $v_0|_{\Omega^c_b}\equiv0$, it follows from $(A_1)$ that
$$
\int_V(\lambda a|u_0|^p+\lambda b|v_0|^q)d\mu=\int_{\Omega_a}\lambda a|u_0|^p+\int_{V\setminus \Omega_a}\lambda a|u_0|^p +\int_{\Omega_b}\lambda b|v_0|^q+\int_{V\setminus \Omega_b}\lambda b|v_0|^q d\mu=0
$$
which shows that $(u_0, v_0)\in W_\lambda$ for all $\lambda> 0$.
Next, we claim that $(u_0,v_0) \neq (0,0)$. In fact, if $(u_0,v_0) = (0,0)$, without loss of generality, setting  $p\geq q$, then by the fact that $(u_{\lambda_k},v_{\lambda_k}) \in \mathcal{N}_{\lambda_k}$ and $F(x,0,0)=0$ for all $x\in V$, (\ref{t7}) and (\ref{a3}), we have
\begin{eqnarray*}
0<M&=&\lim\limits_{k \rightarrow \infty}J_{\lambda_k}(u_{\lambda_k},v_{\lambda_k})\\
&=&\lim\limits_{k \rightarrow \infty}\left[\frac{1}{p}\|u_{\lambda_k}\|^p_{W_\lambda(a)} + \frac{1}{q}\|v_{\lambda_k}\|^q_{W_\lambda(b)}- \int_V F(x,u_{\lambda_k},v_{\lambda_k})d\mu\right]\\
&=&\lim\limits_{k \rightarrow \infty}\left[\frac{1}{q}\|u_{\lambda_k}\|^p_{W_\lambda(a)} + \frac{1}{q}\|v_{\lambda_k}\|^q_{W_\lambda(b)} +\left(\frac{1}{p}-\frac{1}{q}\right)\|u_{\lambda_k}\|^p_{W_\lambda(a)}- \int_V F(x,u_{\lambda_k},v_{\lambda_k})d\mu\right]\\
&=&\lim\limits_{k \rightarrow \infty}\left[\frac{1}{q}\int_V F_{u_{\lambda_k}}(x,u_{\lambda_k},v_{\lambda_k})u_{\lambda_k}d\mu+ \frac{1}{q}\int_V F_{v_{\lambda_k}}(x,u_{\lambda_k},v_{\lambda_k})v_{\lambda_k}d\mu\right.\\
&&+\left.\left(\frac{1}{p}-\frac{1}{q}\right)\|u_{\lambda_k}\|^p_{W_\lambda(a)}- \int_V F(x,u_{\lambda_k},v_{\lambda_k})d\mu\right]\\
&=&\frac{1}{q}\int_V \left(F_{u_{0}}(x,u_{0},v_0)u_{0}+F_{v_0}(x,u_0,v_0)v_0\right)d\mu
+\lim\limits_{k \rightarrow \infty}\left(\frac{1}{p}-\frac{1}{q}\right)\|u_{\lambda_k}\|^p_{W_\lambda(a)}- \int_V F(x,u_0,v_0)d\mu\\
&=&\lim\limits_{k \rightarrow \infty}\left(\frac{1}{p}-\frac{1}{q}\right)\|u_{\lambda_k}\|^p_{W_\lambda(a)}\\
&\leq&0.
\end{eqnarray*}
It is a contradiction.  Hence $(u_0,v_0) \neq (0,0)$.
\par
Combining Lemma 3.3, it is easy to verify that there exists a constant $t'>0$ such that $(t'u_0,t'v_0)\in \mathcal{N}_\lambda$ and further by $u_0|_{\Omega^c_a}\equiv0$ and $v_0|_{\Omega^c_b}\equiv0$, we can obtain $(t'u_0,t'v_0)\in \mathcal{N}_\Omega$. Then by the weak lower semi-continuity of norm, (\ref{t5}) and (\ref{t7}), we can infer that
\begin{eqnarray*}
m_\Omega&\leq&J_\Omega(t'u_0,t'v_0)\\
&=&\frac{1}{p}\left(\int_{\overline{\Omega}_a}|\nabla (tu_0)|^pd\mu+\int_{\Omega_a}|t'u_0|^pd\mu\right)+\frac{1}{q}\left(\int_{\overline{\Omega}_b} |\nabla(t'v_0)|^q d\mu+\int_{\Omega_b}|t'v_0|^q d\mu\right)\nonumber\\
&&-\int_{\Omega_a\cup\Omega_b} F(x,t'u_0,t'v_0)d\mu\\
&=&\frac{1}{p}\int_V (t')^p (|\nabla  u_0|^p+| u_0|^p)d\mu+ \frac{1}{q}\int_V (t')^q (|\nabla  v_0|^q+| v_0|^q)d\mu-\int_V F(x,t'u_0,t'v_0)d\mu\\
&\leq& \liminf\limits_{k \rightarrow \infty} \left[\frac{1}{p}\int_V (t')^p(|\nabla u_{\lambda_k}|^p+ |u_{\lambda_k}|^p)d\mu + \frac{1}{q}\int_V (t')^q(|\nabla v_{\lambda_k}|^q+ |v_{\lambda_k}|^q)d\mu
-\int_V F(x,t'u_{\lambda_k},t'v_{\lambda_k})d\mu\right]\\
&\leq& \liminf\limits_{k \rightarrow \infty} \left[\frac{1}{p}\int_V (t')^p(|\nabla u_{\lambda_k}|^p+ (\lambda_k a+1)|u_{\lambda_k}|^p)d\mu + \frac{1}{q}\int_V (t')^q(|\nabla v_{\lambda_k}|^q+ (\lambda_k b +1)|v_{\lambda_k}|^q)d\mu\right.\\
& &   \left.-\int_V F(x,t'u_{\lambda_k},t'v_{\lambda_k})d\mu\right]\\
&=&\liminf\limits_{k \rightarrow \infty} J_{\lambda_k}(t'u_{\lambda_k}, t'v_{\lambda_k})\\
&\leq& \liminf\limits_{k \rightarrow \infty} J_{\lambda_k}( u_{\lambda_k}, v_{\lambda_k}) \;\;\;\;\mbox{(by Lemma 3.4)}\\
&=&M.
\end{eqnarray*}
Then combining with (\ref{Z10}), we get
$
\lim\limits_{\lambda \rightarrow \infty}m_\lambda=m_\Omega.
$
The proof is complete.
\qed

\vskip2mm
\noindent
{\bf Proof of Theorem 1.2. }
From Lemma 4.1, $(t'u_0,t'v_0)\in \mathcal{N}_\lambda$ and $t'>0$, we obtain that $(u_0,v_0) \not\equiv (0,0)$. We have proved in Lemma 4.1 that $u_0|_{\Omega^c_a}=0$ and $v_0|_{\Omega^c_b}=0$. Then there exists some $t>0$ such that $(t'u_0,t'v_0) \in \mathcal{N}_\Omega$. First we claim that as $k \rightarrow + \infty$, there hold
\begin{eqnarray}
\label{f8}
\lambda_k\int_V a|u_{\lambda_k}|^pd\mu \rightarrow 0,  \;\;\; \lambda_k\int_V b|v_{\lambda_k}|^qd\mu  \rightarrow 0
\end{eqnarray}
and
\begin{eqnarray}
\label{f6}
\int_V (|\nabla u_{\lambda_k}|^p +|u_{\lambda_k}|^p)d\mu  \rightarrow\int_V( |\nabla  u_0|^p+|u_0|^p)d\mu, \;\;\;
 \int_V (|\nabla  v_{\lambda_k}|^q +|v_{\lambda_k}|^p )d\mu \rightarrow \int_V (|\nabla  v_0|^q +|v_0|^q)d\mu.
\end{eqnarray}
Otherwise, by (\ref{t5}) and the weak lower semi-continuity of the norm, for some $\eta_1>0$ and $\eta_2>0$, there holds
$$
 \lim\limits_{k \rightarrow \infty} \left(\lambda_k\int_V a|u_{\lambda_k}|^pd\mu\right)= \eta_1\;\;\;\;\; \mbox{or}\;\;\;\;\;\lim\limits_{k \rightarrow \infty} \left(\lambda_k\int_V b|v_{\lambda_k}|^qd\mu\right) = \eta_2
$$
or
$$
\liminf\limits_{k \rightarrow \infty} \left(\int_V( |\nabla u_{\lambda_k}|^p+|  u_{\lambda_k}|^p)d\mu \right) > \int_V (|\nabla  u_0|^p+ |u_0|^p)d\mu
$$
or
$$
\liminf\limits_{k \rightarrow \infty} \left(\int_V( |\nabla v_{\lambda_k}|^q+|  v_{\lambda_k}|^q)d\mu \right) > \int_V (|\nabla  v_0|^q+ |v_0|^q)d\mu.
$$
If any of the above four formulas hold, we get
\begin{eqnarray*}
&&J_\Omega(t'u_0,t'v_0)\\
&=&\frac{1}{p}\left(\int_{\overline{\Omega}_a} |\nabla (t'u_0)|^pd\mu+\int_{\Omega_a}|t'u_0|^p d\mu\right)+\frac{1}{q}\left(\int_{\overline{\Omega}_b} |\nabla(t'v_0)|^qd\mu+\int_{\Omega_b}|t'v_0|^q d\mu\right)\nonumber\\
&&-\int_{\Omega_a\cup\Omega_b} F(x,t'u_0,t'v_0)d\mu\\
&=&\frac{1}{p}\int_V (t')^p (|\nabla  u_0|^p+| u_0|^p)d\mu+ \frac{1}{q}\int_V (t')^q (|\nabla  v_0|^q+| v_0|^q)d\mu-\int_V F(x,t'u_0,t'v_0)d\mu\\
&<& \liminf\limits_{k \rightarrow \infty} \left[\frac{1}{p}\int_V (t')^p(|\nabla u_{\lambda_k}|^p+ (\lambda_k a+1)|u_{\lambda_k}|^p)d\mu
+ \frac{1}{q}\int_V (t')^q(|\nabla v_{\lambda_k}|^q+ (\lambda_k b +1)|v_{\lambda_k}|^q)d\mu\right.\\
&&\left. -\int_V F(x,t'u_{\lambda_k},t'v_{\lambda_k})d\mu\right]\\
&=&\liminf\limits_{k \rightarrow \infty} J_{\lambda_k}(t'u_{\lambda_k}, t'v_{\lambda_k})\\
&\leq& \liminf\limits_{k \rightarrow \infty} J_{\lambda_k}( u_{\lambda_k}, v_{\lambda_k}) \;\;\;\;\mbox{(by Lemma 3.3)}\\
&=&m_\Omega,
\end{eqnarray*}
which contradicts to the definition of $m_\Omega$ and then the claim is proved.
\par
Now we can prove that $(u_0,v_0)$ is a ground state solution of (\ref{e5}). Firstly, we shall prove $(u_0,v_0)$ is a weak solution of (\ref{e5}). Since $(u_{\lambda_k},v_{\lambda_k})$ is the ground state solution, then $J'_{\lambda_k}(u_{\lambda_k},v_{\lambda_k})=0$. Thus for any $(\phi_1, \phi_2) \in   W_{\Omega_a}\times W_{\Omega_b}\subset W_\lambda $, we have
\begin{eqnarray}
\label{X8}
&&\int_V(|\nabla u_{\lambda_k}|^{p-2}\Gamma(u_{\lambda_k},\phi_1)+(\lambda_k a+1)|u_{\lambda_k}|^{p-2}u_{\lambda_k}\phi_1)d\mu
+\int_V(|\nabla v_{\lambda_k}|^{q-2}\Gamma( v_{\lambda_k}, \phi_2)+(\lambda_k b+1)|v_{\lambda_k}|^{q-2}v_{\lambda_k}\phi_2)d\mu\nonumber\\
&=&\int_V (F_{u_{\lambda_k}}(x,u_{\lambda_k},v_{\lambda_k})\phi_1+F_{v_{\lambda_k}}(x,u_{\lambda_k},v_{\lambda_k})\phi_2)d\mu.
\end{eqnarray}
By $(A_1)$ and noting that $\Omega_a=\{x\in V, a(x)=0\}$ and $\Omega_b=\{x\in V, b(x)=0\}$, it is easy to obtain that $a(x)\phi_1(x) \equiv 0$ and $b(x)\phi_2(x)\equiv 0 $ for any $x \in V$. Then (\ref{X8}) reduces to
\begin{eqnarray}
\label{XX5}
&&\int_{\overline{\Omega}_a} |\nabla u_{\lambda_k}|^{p-2}\Gamma( u_{\lambda_k} ,\phi_1)d\mu+\int_{\Omega_a}|u_{\lambda_k}|^{p-2}u_{\lambda_k}\phi_1 d\mu
+\int_{\overline{\Omega}_b} |\nabla v_{\lambda_k}|^{q-2}\Gamma( v_{\lambda_k},\phi_2)d\mu+\int_{\Omega_b}|v_{\lambda_k}|^{q-2}v_{\lambda_k}\phi_2 d\mu\nonumber\\
&=&\int_{\Omega_a\cup\Omega_b}  F_{u_{\lambda_k}}(x,u_{\lambda_k},v_{\lambda_k})\phi_1d\mu+\int_{\Omega_a\cup\Omega_b}F_{v_{\lambda_k}}(x,u_{\lambda_k},v_{\lambda_k})\phi_2d\mu.
\end{eqnarray}
Note that
\begin{eqnarray}
\label{Z12}
&&\int_{\overline{\Omega}_a} |\nabla u_{\lambda_k}|^{p-2}\Gamma( u_{\lambda_k} ,\phi_1)d\mu-\int_{\overline{\Omega}_a} |\nabla u_0|^{p-2}\Gamma( u_0 ,\phi_1)d\mu\nonumber\\
&=&\int_{\overline{\Omega}_a} (|\nabla u_{\lambda_k}|^{p-2}\Gamma( u_{\lambda_k} ,\phi_1)d\mu-|\nabla u_{\lambda_k}|^{p-2}\Gamma( u_{0} ,\phi_1))d\mu
+\int_{\overline{\Omega}_a} (|\nabla u_{\lambda_k}|^{p-2}\Gamma( u_{0} ,\phi_1)-|\nabla u_{0}|^{p-2}\Gamma( u_{0} ,\phi_1)) d\mu\nonumber\\
&=&\int_{\overline{\Omega}_a} (|\nabla u_{\lambda_k}|^{p-2}\Gamma( u_{\lambda_k}-u_0 ,\phi_1)d\mu +\int_{\overline{\Omega}_a} (|\nabla u_{\lambda_k}|^{p-2}-|\nabla u_{0}|^{p-2})\Gamma( u_{0} ,\phi_1) d\mu.
\end{eqnarray}
By (\ref{Z10}), $\|(u_{\lambda_k},v_{\lambda_k})\|_{W_\Omega} \leq L$. According to the definition of the local finite graph, we obtain that $\Omega_a, \Omega_b, \partial\Omega_a$ and $\partial\Omega_b$ are finite sets. So $W_\Omega$ is a finite dimensional space. Then, up to a subsequence, $(u_{\lambda_k},v_{\lambda_k})\rightarrow (u_0,v_0)$ in $W_\Omega$. Thus we can get $\int_{\overline{\Omega}_a} |\nabla( u_{\lambda_k}-u_0)|^p d\mu\rightarrow 0$ and $\int_{\overline{\Omega}_b} |\nabla( v_{\lambda_k}-v_0)|^ q d\mu\rightarrow 0$. Similar to $(7.10)-(7.11)$ in \cite{Yang P 2023}, we can get $\int_{\overline{\Omega}_a} |\nabla u_{\lambda_k}|^{p-2}\Gamma( u_{\lambda_k}-u_0 ,\phi_1)d\mu \rightarrow0$ and $\int_{\overline{\Omega}_a} (|\nabla u_{\lambda_k}|^{p-2}-|\nabla u_{0}|^{p-2})\Gamma( u_{0} ,\phi_1) d\mu \rightarrow0$. Thus (\ref{Z12}) shows that
\begin{eqnarray*}
\label{X1}
\lim\limits_{k \rightarrow \infty} \int_{\overline{\Omega}_a} |\nabla u_{\lambda_k}|^{p-2}\Gamma( u_{\lambda_k} ,\phi_1)d\mu=\int_{\overline{\Omega}_a} |\nabla u_0|^{p-2}\Gamma( u_0 ,\phi_1)d\mu,
\end{eqnarray*}
and similarly, we also have
\begin{eqnarray*}
\label{XX1}
\lim\limits_{k \rightarrow \infty} \int_{\overline{\Omega}_b} |\nabla v_{\lambda_k}|^{q-2}\Gamma( v_{\lambda_k} ,\phi_2)d\mu=\int_{\overline{\Omega}_b} |\nabla v_0|^{q-2}\Gamma( v_0 ,\phi_2)d\mu.
\end{eqnarray*}
Using (\ref{Z10}) again, we have
$$
L  \geq\|u_{\lambda_k}\|_{W_{\Omega_a}}
\geq \left( \int_{{\Omega}_a} | u_{\lambda_k}|^pd\mu\right)^{\frac{1}{p}}
\geq |{\Omega}_a|\cdot\|u_{\lambda_k}\|_{\infty,{\Omega}_a}^p.
$$
So $|u_{\lambda_k}| \leq \left(\frac{L}{|{\Omega}_a|}\right)^{\frac{1}{p}}$ for all $x\in {\Omega}_a$. Similarly, we can also get $|v_{\lambda_k}| \leq \left(\frac{L}{|{\Omega}_b|}\right)^{\frac{1}{q}}$ for all $x\in {\Omega}_b$.
Note that $\Omega_a, \Omega_b$ and $\Omega_a\cup \Omega_b$ are finite sets. Then by (\ref{t5}) and $(F_1)$, we have
\begin{eqnarray*}
\label{X4}
&&\lim\limits_{k \rightarrow \infty}\int_{\Omega_a} |u_{\lambda_k}|^{p-2}u_{\lambda_k}\phi_1 d\mu=\int_{\Omega_a} |u_0|^{p-2}u_0\phi_1d\mu ,\\
&&\lim\limits_{k \rightarrow \infty}\int_{\Omega_b} |v_{\lambda_k}|^{q-2}v_{\lambda_k}\phi_2 d\mu=\int_{\Omega_b} |v_0|^{q-2}v_0\phi_2d\mu,\\
&&\lim\limits_{k \rightarrow \infty}\int_{\Omega_a\cup\Omega_b} F_{u_{\lambda_k}}(x,u_{\lambda_k},v_{\lambda_k})\phi_1d\mu = \int_{\Omega_a\cup\Omega_b} F_{u_{0}}(x,u_{0},v_{0})\phi_1d\mu,\\
&&\lim\limits_{k \rightarrow \infty}\int_{\Omega_a\cup\Omega_b} F_{v_{\lambda_k}}(x,u_{\lambda_k},v_{\lambda_k})\phi_2d\mu = \int_{\Omega_a\cup\Omega_b} F_{v_{0}}(x,u_{0},v_{0})\phi_2d\mu.
\end{eqnarray*}
Then as $k \rightarrow \infty$, (\ref{XX5}) becomes
\begin{eqnarray*}
&&\int_{\overline{\Omega}_a} |\nabla u_0|^{p-2}\Gamma( u_0 , \phi_1)d\mu+\int_{\Omega_a}|u_0|^{p-2}u_0\phi_1 d\mu
+\int_{\overline{\Omega}_b} |\nabla v_0|^{q-2}\Gamma( v_0, \phi_2)d\mu+\int_{\Omega_b}|v_0|^{q-2}v_0\phi_2 d\mu\\
&=&\int_{\Omega_a\cup\Omega_b}  F_{u_0}(x,u_0,v_0)\phi_1d\mu+\int_{\Omega_a\cup\Omega_b}F_{v_0}(x,u_0,v_0)\phi_2d\mu,
\end{eqnarray*}
which tells us that $J'_\Omega(u_0,v_0)=0$. Hence $(u_0,v_0)\in  \mathcal{N}_\Omega$ is a solution of (\ref{e5}).
\par
Next, we shall prove $(u_0,v_0)$ is a ground state solution of (\ref{e5}). By $u_0|_{\Omega^c_a}=0, v_0|_{\Omega^c_b}=0$, (\ref{t5}), (\ref{t7}), (\ref{f8}) and (\ref{f6}), we have
\begin{eqnarray*}
&&\lim\limits_{k \rightarrow \infty}J_{\lambda_k}(u_{\lambda_k},v_{\lambda_k})\\
&=&\lim\limits_{k \rightarrow \infty}\left[\frac{1}{p}\int_V(|\nabla u_{\lambda_k}|^p+(\lambda_k a+1)|u_{\lambda_k}|^p)d\mu+\frac{1}{q}\int_V(|\nabla v_{\lambda_k}|^q+(\lambda_k b+1)|v_{\lambda_k}|^q)d\mu
-\int_V F(x,u_{\lambda_k},v_{\lambda_k})d\mu\right]\\
&=&\frac{1}{p}\int_V(|\nabla u_0|^p+|u_0|^p)d\mu+\frac{1}{q}\int_V(|\nabla v_0|^q+|v_0|^q)d\mu
-\int_V F(x,u_0,v_0)d\mu+o_k(1)\\
&=&\frac{1}{p}\left(\int_{\overline{\Omega}_a}|\nabla u_0|^pd\mu+\int_{\Omega_a}|u_0|^pd\mu\right)+\frac{1}{q}\left(\int_{\overline{\Omega}_b}|\nabla v_0|^qd\mu+\int_{\Omega_b}|v_0|^qd\mu\right)
-\int_{\Omega_a\cup\Omega_b} F(x,u_0,v_0)d\mu+o_k(1)\\
&=&J_\Omega(u_0,v_0)+o_k(1).
\end{eqnarray*}
Since $\lim\limits_{k \rightarrow \infty}J_{\lambda_k}( u_{\lambda_k}, v_{\lambda_k})=\lim\limits_{k \rightarrow \infty}m_{\lambda_k}$ and Lemma 4.1 tells us that $\lim\limits_{k \rightarrow \infty} m_{\lambda_k}=m_\Omega$, then $J_\Omega(u_0,v_0)=m_\Omega$. Thus we get that $(u_0,v_0)$ is a solution of (\ref{e5}) which achieves the ground state. Finally, by Lemma 4.1 and the above proofs, we can conclude that for any sequence $\lambda_k \rightarrow +\infty$, up to a subsequence, the corresponding ground state solutions $(u_{\lambda_k},v_{\lambda_k}) \in \mathcal{N}_{\lambda_k}$ of (\ref{eq1}) satisfying $J_{\lambda_k}(u_{\lambda_k},v_{\lambda_k}) =m_{\lambda_k}$ converge to a ground state solution of (\ref{e5}) with ground state energy $m_\Omega$. Thus Theorem 1.2 is proved.
\qed

\vskip2mm

{\section{Appendix A}}
  \setcounter{equation}{0}
\noindent
 {\bf Appendix A.1.}  {\it Assume that $\frac{1}{\gamma}<\beta +\frac{1}{\beta} -1$. Then $ \frac{1}{1+\max\{\frac{1}{p}+\frac{q-1}{q},  \frac{1}{q}+\frac{p-1}{p}\}} >\frac{1}{1+\beta}$.}

 \vskip0mm
 \noindent
{\bf Proof.}\ \ When $p\geq q, \beta =p, \gamma =q$, we have
$$
\max\left\{\frac{1}{p}+\frac{q-1}{q},  \frac{1}{q}+\frac{p-1}{p}\right\}=\frac{1}{q}+\frac{p-1}{p},
$$
then
$$
\frac{1}{q}+\frac{p-1}{p}-p=\frac{p+pq-q-p^2 q}{pq}=\frac{p-q+pq(1-p)}{pq}.
$$
If $p-q+pq(1-p) <0$, then $1-p <\frac{1}{p}-\frac{1}{q}$, we can obtain $\frac{1}{q} <p+ \frac{1}{p}-1$.
\par
When $p\leq q, \beta =q, \gamma =p$, we have
$$
\max\{\frac{1}{p}+\frac{q-1}{q},  \frac{1}{q}+\frac{p-1}{p}\}=\frac{1}{p}+\frac{q-1}{q},
$$
then
$$
\frac{1}{p}+\frac{q-1}{q}-q=\frac{q+pq-p-p q^2}{pq}=\frac{q-p+pq(1-q)}{pq}.
$$
If $q-p+pq(1-q) <0$, then $1-q <\frac{1}{q}-\frac{1}{p}$, we can obtain $\frac{1}{p} <q+ \frac{1}{q}-1$.\\
Thus when $\frac{1}{\gamma}<\beta +\frac{1}{\beta} -1$, we have $\max\{\frac{1}{p}+\frac{q-1}{q}, \frac{1}{q}+\frac{p-1}{p}\}<\beta $ i.e
 $ \frac{1}{1+\max\{\frac{1}{p}+\frac{q-1}{q},  \frac{1}{q}+\frac{p-1}{p}\}} >\frac{1}{1+\beta}$.
\qed

 \vskip2mm
\noindent
{\bf Appendix A.2.} {\it Assume that $(F_1)$ and $(F_3)$ hold. Then $$|F(x,t,s)|\leq \frac{2C_1(x)}{p}|t|^p+\left(\frac{C_1(x)(p-1)}{p}+\frac{C_1(x)}{q}\right)|s|^q+\frac{2C_2(x)}{r_1}|t|^{r_1}+\left(\frac{C_2(x)(r_1-1)}{r_1}+\frac{C_2(x)}{r_2}\right)|s|^{r_2} $$
for all $x\in V$ and $\forall (t,s) \in \R^2$.}
\vskip2mm
\noindent
{\bf Proof.}\ \ From $(F_3)$ and Young inequality, we can get
\begin{eqnarray*}
&&|F(x,t,s)|\\
&\leq&  \int^{|t|}_0|F_t(x,t,s)|dt+\int^{|s|}_0|F_s(x,0,s)|ds\\
&\leq&  \int^{|t|}_0[ C_1(x)(|t|^{p-1}+|s|^{\frac{q(p-1)}{p}})+ C_2(x)(|t|^{r_1-1}+|s|^{\frac{r_2(r_1-1)}{r_1}})]dt+\int^{|s|}_0 (C_1(x)|s|^{q-1}+C_2(x)|s|^{r_2-1})ds\\
&=  &   \frac{C_1(x)}{p}|t|^{p}+ C_1(x)|s|^{\frac{q(p-1)}{p}}|t| +  \frac{C_2(x)}{r_1}|t|^{r_1}+ C_2(x)|s|^{\frac{r_2(r_1-1)}{r_1}}|t| +\frac{C_1(x)}{q} |s|^{q}+\frac{C_2(x)}{r_2}|s|^{r_2}\\
&\leq&  \frac{C_1(x)}{p}|t|^{p}+ \frac{C_1(x)(p-1)}{p} |s|^q +\frac{C_1(x)}{p}|t|^{p} +\frac{C_2(x)}{r_1}|t|^{r_1}
+ \frac{C_2(x)(r_1-1)}{r_1}|s|^{r_2} + \frac{C_2(x)}{r_1}|t|^{r_1}
+\frac{C_1(x)}{q} |s|^{q}+\frac{C_2(x)}{r_2}|s|^{r_2}\\
&=  &\frac{2C_1(x)}{p}|t|^p+\left(\frac{C_1(x)(p-1)}{p}+\frac{C_1(x)}{q}\right)|s|^q+\frac{2C_2(x)}{r_1}|t|^{r_1}+\left(\frac{C_2(x)(r_1-1)}{r_1}+\frac{C_2(x)}{r_2}\right)|s|^{r_2}.
\end{eqnarray*}
Thus, the proof is finished.\qed

 \vskip2mm
\noindent
{\bf Appendix A.3.} {\it Assume that $(F_3)$ holds. Then
$$
|F_{t}(x,t,s)t+F_{s}(x,t,s)s|  \leq  C_3(x)(|t|^p+|s|^q)+C_4(x)(|t|^{r_1}+|s|^{r_2})
$$
for all $x\in V$ and $\forall (t,s) \in \R^2$, where $C_3(x)=C_1(x)\cdot\max\left\{1+\frac{1}{p}+\frac{q-1}{q}, 1+\frac{1}{q}+\frac{p-1}{p}\right\} $ and $C_4(x)=C_2(x)\cdot\max\left\{1+\frac{1}{r_1}+\frac{r_2-1}{r_2}, 1+\frac{1}{r_2}+\frac{r_1-1}{r_1}\right\}$.}
\vskip2mm
\noindent
{\bf Proof.} \ It follows from  $(F_3)$  that
\begin{eqnarray*}
&&|F_{t}(x,t,s)t+F_{s}(x,t,s)s| \\
&\leq&  C_1(x)(|t|^{p}+|s|^{\frac{q(p-1)}{p}}|t|)+C_2(x)(|t|^{r_1}+|s|^{\frac{r_2(r_1-1)}{r_1}}|t|)
+ C_1(x)(|t|^{\frac{p(q-1)}{q}}|s|+|s|^{q})+C_2(x)(|t|^{\frac{r_1(r_2-1)}{r_2}}|s|+|s|^{r_2})\\
&\leq&  C_1(x)|t|^{p} +\frac{C_1(x)(p-1)}{p}|s|^{q}+\frac{C_1(x)}{p}|t|^{p} + C_2(x)|t|^{r_1}+  \frac{C_2(x)(r_1-1)}{r_1}|s|^{r_2}+\frac{C_2(x)}{r_1}|t|^{r_1}\\
&&+\frac{C_1(x)(q-1)}{q}|t|^{p}+ \frac{C_1(x)}{q}|s|^{q}  +C_1(x)|s|^{q} + \frac{C_2(x)(r_2-1)}{r_2}|t|^{r_1}+\frac{C_2(x)}{r_2}|s|^{r_2} + C_2(x)|s|^{r_2}\\
&=  & \left(C_1(x)+\frac{C_1(x)}{p}+\frac{C_1(x)(q-1)}{q}\right)|t|^p  +  \left(C_1(x)+\frac{C_1(x)}{q}+\frac{C_1(x)(p-1)}{p}\right)|s|^q \\
&& + \left(C_2(x)+\frac{C_2(x)}{r_1}+\frac{C_2(x)(r_2-1)}{r_2}\right)|t|^{r_1}  +  \left(C_2(x)+\frac{C_2(x)}{r_2}+\frac{C_2(x)(r_1-1)}{r_1}\right)|s|^{r_2} \\
&\leq&   C_1(x)\cdot\max\left\{1+\frac{1}{p}+\frac{q-1}{q}, 1+\frac{1}{q}+\frac{p-1}{p}\right\}(|t|^p +|s|^q ) \\
&&+C_2(x)\cdot\max\left\{1+\frac{1}{r_1}+\frac{r_2-1}{r_2}, 1+\frac{1}{r_2}+\frac{r_1-1}{r_1}\right\}(|t|^{r_1}+|s|^{r_2})\\
&=  & C_3(x)(|t|^p+|s|^q)+C_4(x)(|t|^{r_1}+|s|^{r_2}).
\end{eqnarray*}
Thus, we finish the proof. \qed

 \vskip2mm
\noindent
{\bf Appendix A.4.} {\it If $\beta=p, t_\lambda \in (0,1)$ and $(F_4)$ holds, then (\ref{t9}) holds.}
\vskip2mm
\noindent
{\bf Proof.}\ \ Let $A(t_\lambda)=\frac{1}{p} F_{t_\lambda u_\lambda}(x, t_\lambda u_\lambda ,t_\lambda v_\lambda )t_\lambda u_\lambda +\frac{1}{p} F_{t_\lambda v_\lambda}(x,t_\lambda u_\lambda ,t_\lambda v_\lambda )t_\lambda v_\lambda - F(x,t_\lambda u_\lambda ,t_\lambda v_\lambda )$.\\
Then
\begin{eqnarray*}
&&A'(t_\lambda)\\
& =   & \frac{1}{p}F_{t_\lambda u_\lambda,t_\lambda u_\lambda}(x, t_\lambda u_\lambda ,t_\lambda v_\lambda )t_\lambda u_\lambda^2+\frac{1}{p}F_{t_\lambda u_\lambda,t_\lambda v_\lambda}(x, t_\lambda u_\lambda ,t_\lambda v_\lambda )t_\lambda u_\lambda v_\lambda+\frac{1}{p} F_{t_\lambda u_\lambda}(x, t_\lambda u_\lambda ,t_\lambda v_\lambda )u_\lambda \\
&&+\frac{1}{p}F_{t_\lambda v_\lambda,t_\lambda u_\lambda}(x, t_\lambda u_\lambda ,t_\lambda v_\lambda )t_\lambda u_\lambda v_\lambda+\frac{1}{p}F_{t_\lambda v_\lambda,t_\lambda v_\lambda}(x, t_\lambda u_\lambda ,t_\lambda v_\lambda )t_\lambda v_\lambda^2+\frac{1}{p} F_{t_\lambda v_\lambda}(x, t_\lambda u_\lambda ,t_\lambda v_\lambda )v_\lambda\\
&&-F_{t_\lambda u_\lambda}(x, t_\lambda u_\lambda ,t_\lambda v_\lambda ) u_\lambda-F_{t_\lambda v_\lambda}(x,t_\lambda u_\lambda ,t_\lambda v_\lambda )v_\lambda\\
& >   & \frac{p-1}{p}F_{t_\lambda u_\lambda}(x, t_\lambda u_\lambda ,t_\lambda v_\lambda ) u_\lambda +\frac{1}{p}F_{t_\lambda u_\lambda}(x, t_\lambda u_\lambda ,t_\lambda v_\lambda ) u_\lambda +\frac{1}{p}F_{t_\lambda u_\lambda,t_\lambda v_\lambda}(x, t_\lambda u_\lambda ,t_\lambda v_\lambda )t_\lambda u_\lambda v_\lambda\\
&&+\frac{p-1}{p}F_{t_\lambda v_\lambda}(x, t_\lambda u_\lambda ,t_\lambda v_\lambda ) v_\lambda+\frac{1}{p} F_{t_\lambda v_\lambda}(x, t_\lambda u_\lambda ,t_\lambda v_\lambda )v_\lambda +\frac{1}{p}F_{t_\lambda v_\lambda,t_\lambda u_\lambda}(x, t_\lambda u_\lambda ,t_\lambda v_\lambda )t_\lambda u_\lambda v_\lambda\\
&&-F_{t_\lambda u_\lambda}(x, t_\lambda u_\lambda ,t_\lambda v_\lambda ) u_\lambda-F_{t_\lambda v_\lambda}(x,t_\lambda u_\lambda ,t_\lambda v_\lambda )v_\lambda\\
& =   & \frac{2}{p}F_{t_\lambda u_\lambda,t_\lambda v_\lambda}(x, t_\lambda u_\lambda ,t_\lambda v_\lambda )t_\lambda u_\lambda v_\lambda\\
& \geq   & 0.
\end{eqnarray*}
So $A(t_\lambda)$ is a strictly monotonically increasing function. Thus $A(1)>A(t_\lambda)$.
\qed

 \vskip2mm
\noindent
{\bf Appendix A.5.} {\it $\mathcal{N}_\Omega  \subseteq  \mathcal{N}_\lambda $}.
\vskip2mm
\noindent
{\bf Proof.}\ \ For any $(u_0,v_0)\in\mathcal{N}_{\Omega}$, we have
\begin{eqnarray*}
\label{Z11}
\int_{\overline{\Omega}_a}|\nabla u_0|^{p} d\mu+\int_{\Omega_a}|u_0|^{p}d\mu
+\int_{\overline{\Omega}_b}|\nabla v_0|^{q} d\mu+\int_{\Omega_b}|v_0|^{q}d\mu
-\int_{\Omega_a\cup\Omega_b} F_{u_0}(x,u_0,v_0)u_0d\mu-\int_{\Omega_a\cup\Omega_b} F_{v_0}(x,u_0,v_0)v_0d\mu
=0.
\end{eqnarray*}
Then by $(F_1)$, $(A_1)$ and the fact that $W_{\Omega_a}(W_{\Omega_b})$ is the completion of $C_c(\Omega_a)(C_c(\Omega_b))$, we have
\begin{eqnarray*}
\label{e3}
&&\langle J'_\lambda(u_0,v_0),(u_0,v_0)\rangle\nonumber\\
&=&\int_V(|\nabla u_0|^{p } +(\lambda a+1)|u_0|^{p } )d\mu
+\int_V(|\nabla v_0|^{q }  +(\lambda b+1)|v_0|^{q } )d\mu
-\int_V (F_{u_0}(x,u_0,v_0)u_0+F_{v_0}(x,u_0,v_0)v_0)d\mu\nonumber\\
&=&\int_{V \backslash {\overline{\Omega}_a}} |\nabla u_0|^{p }d\mu+\int_{{\Omega}_a} |\nabla u_0|^{p }d\mu+\int_{\partial{\Omega}_a} |\nabla u_0|^{p }d\mu
+ \int_{V \backslash {{\Omega}_a}} (\lambda a+1)| u_0|^{p }d\mu+\int_{{\Omega}_a}  (\lambda a+1)| u_0|^{p }d\mu\nonumber\\
&&+\int_{V \backslash {\overline{\Omega}_b}} |\nabla v_0|^{q }d\mu+\int_{{\Omega}_b} |\nabla v_0|^{q }d\mu+\int_{\partial{\Omega}_b} |\nabla v_0|^{q }d\mu
+ \int_{V \backslash {{\Omega}_b}} (\lambda b+1)| v_0|^{q }d\mu+\int_{{\Omega}_b}  (\lambda b+1)| v_0|^{q }d\mu\nonumber\\
&&-\int_{V \backslash {(\Omega_a \cup {\Omega_b}})}F_{u_0}(x,u_0,v_0)u_0d\mu - \int_{ {\Omega_a \cup {\Omega_b}}}F_{u_0}(x,u_0,v_0)u_0d\mu\nonumber\\
&&-\int_{V \backslash {(\Omega_a \cup {\Omega_b}})}F_{v_0}(x,u_0,v_0)v_0d\mu - \int_{ {\Omega_a \cup {\Omega_b}}}F_{v_0}(x,u_0,v_0)v_0d\mu\nonumber\\
&=&\int_{\overline{\Omega}_a}|\nabla u_0|^{p} d\mu+\int_{\Omega_a}|u_0|^{p}d\mu
+\int_{\overline{\Omega}_b}|\nabla v_0|^{q} d\mu+\int_{\Omega_b}|v_0|^{q}d\mu
-\int_{\Omega_a\cup\Omega_b} (F_{u_0}(x,u_0,v_0)u_0d\mu+ F_{v_0}(x,u_0,v_0)v_0)d\mu\nonumber\\
&=&0.
\end{eqnarray*}
Thus $(u_0,v_0)\in\mathcal{N}_{\lambda}$.
\qed

 \vskip2mm
 \noindent
 {\bf Funding}\\
This project is supported by Yunnan Fundamental Research Projects (grant No: 202301AT070465) and  Xingdian Talent Support Program for Young Talents of Yunnan Province.

\vskip2mm
 \noindent
 {\bf Authors' contributions}\\
The authors contribute the manuscript equally.

 \vskip2mm
 \noindent
 {\bf Competing interests}\\
The authors declare that they have no competing interests.

\vskip2mm
\renewcommand\refname{References}
{}
\end{document}